\documentclass[10pt]{article}
\def\newpic#1{%
   \def\emline##1##2##3##4##5##6{%
         \put(##1,##2){\special{em:point #1##3}}%
         \put(##4,##5){\special{em:point #1##6}}%
         \special{em:line #1##3,#1##6}}}
\newpic{}
\def\emline#1#2#3#4#5#6{%
          \put(#1,#2){\special{em:moveto}}%
          \put(#4,#5){\special{em:lineto}}}
\def\newpic#1{}

\title{Star graphs: threaded distance trees and E-sets}
\author{Italo J. Dejter
\\
{\small Department of Mathematics}
\\
{\small University of Puerto Rico, Rio Piedras, PR 00931-3355}
}
\date{}

\newcommand{\Z}{\hbox{\bf Z}}
\newtheorem{thm}{Theorem}
\newtheorem{propo}[thm]{Proposition}

\newcommand{\subgp}[1]{\langle{#1}\rangle}
\newcommand{\beeq}{\begin{eqnarray*}}
\newcommand{\eneq}{\end{eqnarray*}}
\newcommand{\proof}{\noindent{\bf Proof.\hspace{4mm}}}
\newcommand{\example}{\noindent{\bf Example.\hspace{4mm}}}
\newcommand{\examples}{\noindent{\bf Examples.\hspace{4mm}}}
\newcommand{\rema}{\noindent{\bf Remark.\hspace{4mm}}}
\newcommand{\qed}{\hfill\ \rule{2mm}{2mm}}
\newcommand{\qfd}{\hfill $\fbox{}$\vspace{4mm}}
\begin{document}
\maketitle

\begin{abstract}
The distribution of distances in the star graph $ST_n$, ($1<n\in\Z$), is established, and subsequently a threaded binary tree is obtained that realizes an orientation of $ST_n$ whose levels are given by the distances to the identity permutation, via a pruning algorithm followed by a threading algorithm. In the process, the distributions of distances of the efficient dominating sets of $ST_n$ are determined.
\end{abstract}

\section{Introduction}

The {\it star graph} $ST_n$, ($1<n\in\Z$), is the Cayley graph of the symmetric
group $S_n$ with set of generators $\Theta_n=\{(1\;i),\;
i=2,\ldots,n\}$, (\cite{akers,arum1}). The {\it weight} of a
vertex $u$ of $ST_n$ is its distance to the identity-permutation vertex
$12\ldots n$. In this work, based on  DIMACS Technical Report 2001-05, the weight distributions of certain subsets $C$ of
$ST_n$ are determined, including that of $ST_n$ itself. Theorems 8 and 6 below attain these objectives.
(A variation of Theorem 6 was obtained in a different fashion in \cite{WSLS}).

An independent set $C$ of vertices in a graph is an efficient dominating set \cite{haynes},
or E-set \cite{io}, or 1-perfect codes \cite{Kr}, if each vertex not in $C$ is adjacent to exactly one vertex of $C$. In Section 5, we determine the weight distributions of these E-sets; see Theorem 8 and subsequent remark. In obtaining this, we use a binary directed tree $\Lambda_n=\Lambda(ST_n)$ whose arcs are of two types: {\bf(1)} horizontal, left-to-right, arcs; {\bf(2)} vertical, top-to-bottom, arcs, (as in the subsequent figures). In Section 6, we extend $\Lambda_n$ to an orientation $\Gamma_n$ of $ST_n$, (that is: an oriented graph $\Gamma_n$). Moreover, the graphs $\Gamma_n$ form a nested sequence that converges to a
universal graph $\Gamma_\infty$ associated to the infinite star graph $ST_\infty$.

\section{Definition and examples of $\Lambda_n$}

Let $n>1$ and let $\Sigma\in S_n$. We write $\Sigma=\sigma_1\sigma_2\ldots\sigma_n$, where $\Sigma(i)=\sigma_i$, for $i=1,2,\ldots,n$. A cycle $(\sigma_{i_1}\sigma_{i_2}\ldots\sigma_{i_r})$ of the permutation $\Sigma$
is given by $\Sigma(\sigma_{i_j})=\sigma_{i_{j+1}}$, for $j=1\ldots r$, where $j+1$ is taken as $1$ if $j=r$. Then, $\Sigma$ has {\it length} $r$. Now, $\Sigma$ is said to be {\it proper} if $r>1$. The {\it cycle structure} $\Pi(\Sigma)$ of $\Sigma=\sigma_1\sigma_2\ldots\sigma_n$
is defined as the set of proper cycles of by $\Sigma$.

\newpage

Two vertices $\Sigma^1$ and $\Sigma^2$ of $ST_n$, with 1 in cycles $\tau^1$ of $\Sigma^1$ and $\tau^2$ of $\Sigma^2$ of the same length, have a common {\it $1$-invariant cycle structure} if there is $\Phi\in S_n$ with $\Phi(\Sigma^1)=\Sigma^2$ inducing a 1-1 correspondence $\Phi^*:\Pi(\Sigma^1)\rightarrow\Pi(\Sigma^2)$ sending  $\tau^1$ onto $\tau^2$ and with each $\tau\in\Pi(\Sigma^1)$ and $\Phi(\tau)\in\Pi(\Sigma^2)$ having the same length.
We say that $\Sigma^2$ has the {\it $1$-invariant cycle structure},
(or 1-ics), of $\Sigma^1$. Each vertex $u$ of $\Lambda_n$ is written $$\begin{array}{|c|}\hline w(u),c(u) \\ \Sigma(u) \\ \hline
\end{array}$$ where
{\bf (a)} $\Sigma(u)=\sigma_1\ldots\sigma_{i-1}$ is shorthand for a permutation $\sigma_1\sigma_2\ldots\sigma_n$ of $12\ldots n$ having $i$ as the smallest index in $\{2,\ldots,n\}$ satisfying $\sigma_j=j$, for $i\leq j\leq n$, and $\sigma_j\neq j$ for $1<j<i$;
{\bf (b)} $w(u)$ is the weight of $\Sigma(u)$;
{\bf (c)} $c(u)$ is the cardinality of the set $S(u)$ of permutations having
the 1-ics $\Pi(u)$ of $\Sigma(u)$.

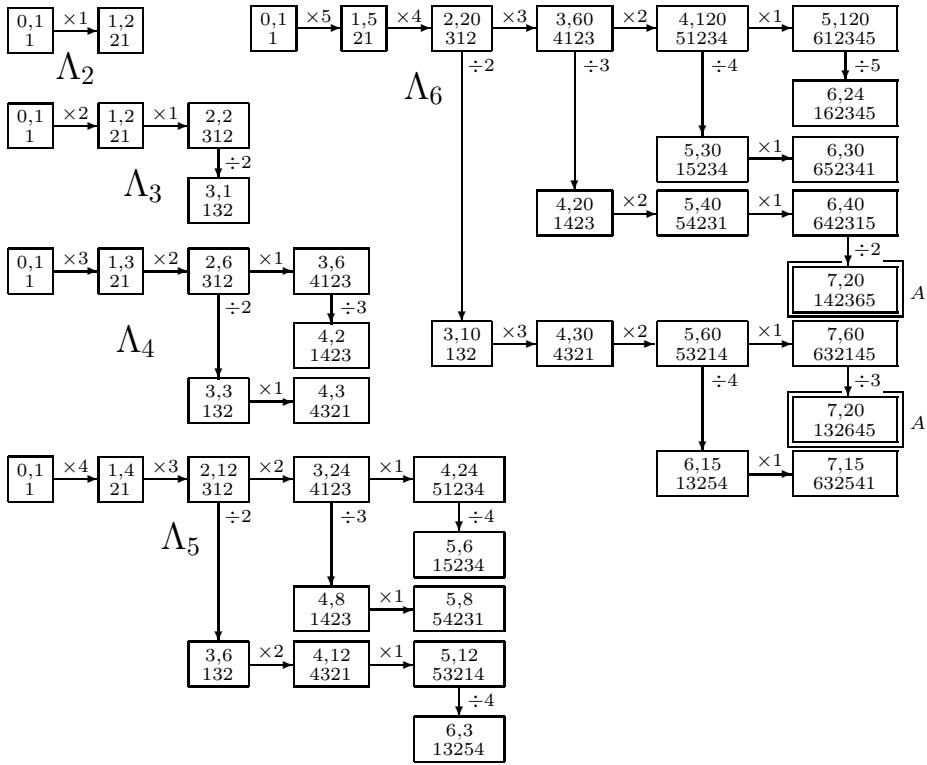
\begin{figure}[htp]
\unitlength=1.00mm
\special{em:linewidth 0.4pt}
\linethickness{0.4pt}
\begin{picture}(121.06,100.33)
\put(3.00,83.03){\makebox(0,0)[cc]{$_{1}$}}
\put(15.00,83.03){\makebox(0,0)[cc]{$_{21}$}}
\put(9.00,85.92){\makebox(0,0)[cc]{$_{\times 2}$}}
\put(28.00,83.03){\makebox(0,0)[cc]{$_{312}$}}
\put(21.00,85.92){\makebox(0,0)[cc]{$_{\times 1}$}}
\put(0.00,87.36){\line(1,0){6.00}}
\put(12.00,87.37){\line(1,0){6.00}}
\put(24.00,87.37){\line(1,0){8.00}}
\put(0.00,81.36){\line(1,0){6.00}}
\put(12.00,81.37){\line(1,0){6.00}}
\put(24.00,81.37){\line(1,0){8.00}}
\put(0.00,87.36){\line(0,-1){6.00}}
\put(6.00,87.37){\line(0,-1){6.00}}
\put(12.00,87.37){\line(0,-1){6.00}}
\put(18.00,87.37){\line(0,-1){6.00}}
\put(24.00,87.37){\line(0,-1){6.00}}
\put(32.00,87.37){\line(0,-1){6.00}}
\put(28.00,73.11){\makebox(0,0)[cc]{$_{132}$}}
\put(24.00,77.45){\line(1,0){8.00}}
\put(24.00,71.45){\line(1,0){8.00}}
\put(24.00,77.45){\line(0,-1){6.00}}
\put(32.00,77.45){\line(0,-1){6.00}}
\put(30.67,79.37){\makebox(0,0)[cc]{$_{\div 2}$}}
\put(3.11,85.37){\makebox(0,0)[cc]{$_{0,1}$}}
\put(15.11,85.37){\makebox(0,0)[cc]{$_{1,2}$}}
\put(28.11,85.37){\makebox(0,0)[cc]{$_{2,2}$}}
\put(28.11,75.45){\makebox(0,0)[cc]{$_{3,1}$}}
\put(3.00,63.69){\makebox(0,0)[cc]{$_{1}$}}
\put(15.00,63.69){\makebox(0,0)[cc]{$_{21}$}}
\put(9.00,66.58){\makebox(0,0)[cc]{$_{\times 3}$}}
\put(28.00,63.69){\makebox(0,0)[cc]{$_{312}$}}
\put(21.00,66.58){\makebox(0,0)[cc]{$_{\times 2}$}}
\put(12.00,68.03){\line(1,0){6.00}}
\put(24.00,68.03){\line(1,0){8.00}}
\put(12.00,62.03){\line(1,0){6.00}}
\put(24.00,62.03){\line(1,0){8.00}}
\put(6.00,68.03){\line(0,-1){6.00}}
\put(12.00,68.03){\line(0,-1){6.00}}
\put(18.00,68.03){\line(0,-1){6.00}}
\put(24.00,68.03){\line(0,-1){6.00}}
\put(32.00,68.03){\line(0,-1){6.00}}
\put(43.00,63.69){\makebox(0,0)[cc]{$_{4123}$}}
\put(35.00,66.58){\makebox(0,0)[cc]{$_{\times 1}$}}
\put(38.00,68.03){\line(0,-1){6.00}}
\put(48.00,68.03){\line(0,-1){6.00}}
\put(38.00,62.03){\line(1,0){10.00}}
\put(38.00,68.03){\line(1,0){10.00}}
\put(43.00,53.79){\makebox(0,0)[cc]{$_{1423}$}}
\put(38.00,58.13){\line(0,-1){6.00}}
\put(48.00,58.13){\line(0,-1){6.00}}
\put(38.00,52.13){\line(1,0){10.00}}
\put(38.00,58.13){\line(1,0){10.00}}
\put(46.00,60.03){\makebox(0,0)[cc]{$_{\div 3}$}}
\put(28.00,46.44){\makebox(0,0)[cc]{$_{132}$}}
\put(24.00,50.78){\line(1,0){8.00}}
\put(24.00,44.78){\line(1,0){8.00}}
\put(24.00,50.78){\line(0,-1){6.00}}
\put(32.00,50.78){\line(0,-1){6.00}}
\put(43.00,46.44){\makebox(0,0)[cc]{$_{4321}$}}
\put(35.00,49.34){\makebox(0,0)[cc]{$_{\times 1}$}}
\put(38.00,50.78){\line(0,-1){6.00}}
\put(48.00,50.78){\line(0,-1){6.00}}
\put(38.00,44.78){\line(1,0){10.00}}
\put(38.00,50.78){\line(1,0){10.00}}
\put(30.67,60.03){\makebox(0,0)[cc]{$_{\div 2}$}}
\put(3.11,66.03){\makebox(0,0)[cc]{$_{0,1}$}}
\put(15.11,66.03){\makebox(0,0)[cc]{$_{1,3}$}}
\put(28.11,66.03){\makebox(0,0)[cc]{$_{2,6}$}}
\put(43.11,66.03){\makebox(0,0)[cc]{$_{3,6}$}}
\put(43.11,56.13){\makebox(0,0)[cc]{$_{4,2}$}}
\put(28.11,48.78){\makebox(0,0)[cc]{$_{3,3}$}}
\put(43.11,48.78){\makebox(0,0)[cc]{$_{4,3}$}}
\put(0.00,68.03){\line(1,0){6.00}}
\put(0.00,62.03){\line(1,0){6.00}}
\put(0.00,68.03){\line(0,-1){6.00}}
\put(3.00,35.99){\makebox(0,0)[cc]{$_{1}$}}
\put(15.00,35.99){\makebox(0,0)[cc]{$_{21}$}}
\put(9.00,38.88){\makebox(0,0)[cc]{$_{\times 4}$}}
\put(28.00,35.99){\makebox(0,0)[cc]{$_{312}$}}
\put(21.00,38.88){\makebox(0,0)[cc]{$_{\times 3}$}}
\put(0.00,40.34){\line(1,0){6.00}}
\put(12.00,40.33){\line(1,0){6.00}}
\put(24.00,40.33){\line(1,0){8.00}}
\put(0.00,34.34){\line(1,0){6.00}}
\put(12.00,34.33){\line(1,0){6.00}}
\put(24.00,34.33){\line(1,0){8.00}}
\put(0.00,40.34){\line(0,-1){6.00}}
\put(6.00,40.33){\line(0,-1){6.00}}
\put(12.00,40.33){\line(0,-1){6.00}}
\put(18.00,40.33){\line(0,-1){6.00}}
\put(24.00,40.33){\line(0,-1){6.00}}
\put(32.00,40.33){\line(0,-1){6.00}}
\put(43.00,35.99){\makebox(0,0)[cc]{$_{4123}$}}
\put(35.00,38.88){\makebox(0,0)[cc]{$_{\times 2}$}}
\put(38.00,40.33){\line(0,-1){6.00}}
\put(48.00,40.33){\line(0,-1){6.00}}
\put(38.00,34.33){\line(1,0){10.00}}
\put(38.00,40.33){\line(1,0){10.00}}
\put(60.00,35.99){\makebox(0,0)[cc]{$_{51234}$}}
\put(51.00,38.88){\makebox(0,0)[cc]{$_{\times 1}$}}
\put(54.00,40.33){\line(0,-1){6.00}}
\put(66.00,40.33){\line(0,-1){6.00}}
\put(54.00,40.33){\line(1,0){12.00}}
\put(54.00,34.33){\line(1,0){12.00}}
\put(60.00,25.93){\makebox(0,0)[cc]{$_{15234}$}}
\put(54.00,30.26){\line(0,-1){6.00}}
\put(66.00,30.26){\line(0,-1){6.00}}
\put(54.00,30.26){\line(1,0){12.00}}
\put(54.00,24.26){\line(1,0){12.00}}
\put(63.00,32.33){\makebox(0,0)[cc]{$_{\div 4}$}}
\put(43.00,18.76){\makebox(0,0)[cc]{$_{1423}$}}
\put(38.00,23.10){\line(0,-1){6.00}}
\put(48.00,23.10){\line(0,-1){6.00}}
\put(38.00,17.10){\line(1,0){10.00}}
\put(38.00,23.10){\line(1,0){10.00}}
\put(60.00,18.76){\makebox(0,0)[cc]{$_{54231}$}}
\put(51.00,21.66){\makebox(0,0)[cc]{$_{\times 1}$}}
\put(54.00,23.10){\line(0,-1){6.00}}
\put(66.00,23.10){\line(0,-1){6.00}}
\put(54.00,23.10){\line(1,0){12.00}}
\put(54.00,17.10){\line(1,0){12.00}}
\put(46.00,32.33){\makebox(0,0)[cc]{$_{\div 3}$}}
\put(28.00,11.41){\makebox(0,0)[cc]{$_{132}$}}
\put(24.00,15.75){\line(1,0){8.00}}
\put(24.00,9.75){\line(1,0){8.00}}
\put(24.00,15.75){\line(0,-1){6.00}}
\put(32.00,15.75){\line(0,-1){6.00}}
\put(43.00,11.41){\makebox(0,0)[cc]{$_{4321}$}}
\put(35.00,14.31){\makebox(0,0)[cc]{$_{\times 2}$}}
\put(38.00,15.75){\line(0,-1){6.00}}
\put(48.00,15.75){\line(0,-1){6.00}}
\put(38.00,9.75){\line(1,0){10.00}}
\put(38.00,15.75){\line(1,0){10.00}}
\put(60.00,11.41){\makebox(0,0)[cc]{$_{53214}$}}
\put(51.00,14.31){\makebox(0,0)[cc]{$_{\times 1}$}}
\put(54.00,15.75){\line(0,-1){6.00}}
\put(66.00,15.75){\line(0,-1){6.00}}
\put(54.00,15.75){\line(1,0){12.00}}
\put(54.00,9.75){\line(1,0){12.00}}
\put(60.00,1.44){\makebox(0,0)[cc]{$_{13254}$}}
\put(54.00,5.78){\line(0,-1){6.00}}
\put(66.00,5.78){\line(0,-1){6.00}}
\put(54.00,5.78){\line(1,0){12.00}}
\put(54.00,-0.22){\line(1,0){12.00}}
\put(63.00,7.75){\makebox(0,0)[cc]{$_{\div 4}$}}
\put(30.67,32.33){\makebox(0,0)[cc]{$_{\div 2}$}}
\put(3.11,38.33){\makebox(0,0)[cc]{$_{0,1}$}}
\put(15.11,38.33){\makebox(0,0)[cc]{$_{1,4}$}}
\put(28.11,38.33){\makebox(0,0)[cc]{$_{2,12}$}}
\put(43.11,38.33){\makebox(0,0)[cc]{$_{3,24}$}}
\put(60.11,38.33){\makebox(0,0)[cc]{$_{4,24}$}}
\put(60.11,28.26){\makebox(0,0)[cc]{$_{5,6}$}}
\put(43.11,21.10){\makebox(0,0)[cc]{$_{4,8}$}}
\put(60.11,21.10){\makebox(0,0)[cc]{$_{5,8}$}}
\put(28.11,13.75){\makebox(0,0)[cc]{$_{3,6}$}}
\put(43.11,13.75){\makebox(0,0)[cc]{$_{4,12}$}}
\put(60.11,13.75){\makebox(0,0)[cc]{$_{5,12}$}}
\put(60.11,3.78){\makebox(0,0)[cc]{$_{6,3}$}}
\put(35.33,95.99){\makebox(0,0)[cc]{$_{1}$}}
\put(47.33,95.99){\makebox(0,0)[cc]{$_{21}$}}
\put(41.33,98.88){\makebox(0,0)[cc]{$_{\times 5}$}}
\put(60.33,95.99){\makebox(0,0)[cc]{$_{312}$}}
\put(53.33,98.88){\makebox(0,0)[cc]{$_{\times 4}$}}
\put(32.33,100.33){\line(1,0){6.00}}
\put(44.33,100.33){\line(1,0){6.00}}
\put(56.33,100.33){\line(1,0){8.00}}
\put(32.33,94.33){\line(1,0){6.00}}
\put(44.33,94.33){\line(1,0){6.00}}
\put(56.33,94.33){\line(1,0){8.00}}
\put(32.33,100.33){\line(0,-1){6.00}}
\put(38.33,100.33){\line(0,-1){6.00}}
\put(44.33,100.33){\line(0,-1){6.00}}
\put(50.33,100.33){\line(0,-1){6.00}}
\put(56.33,100.33){\line(0,-1){6.00}}
\put(64.33,100.33){\line(0,-1){6.00}}
\put(75.33,95.99){\makebox(0,0)[cc]{$_{4123}$}}
\put(67.33,98.88){\makebox(0,0)[cc]{$_{\times 3}$}}
\put(70.33,100.33){\line(0,-1){6.00}}
\put(80.33,100.33){\line(0,-1){6.00}}
\put(70.33,94.33){\line(1,0){10.00}}
\put(70.33,100.33){\line(1,0){10.00}}
\put(92.33,95.99){\makebox(0,0)[cc]{$_{51234}$}}
\put(83.33,98.88){\makebox(0,0)[cc]{$_{\times 2}$}}
\put(86.33,100.33){\line(0,-1){6.00}}
\put(98.33,100.33){\line(0,-1){6.00}}
\put(86.33,100.33){\line(1,0){12.00}}
\put(86.33,94.33){\line(1,0){12.00}}
\put(111.33,95.99){\makebox(0,0)[cc]{$_{612345}$}}
\put(101.33,98.88){\makebox(0,0)[cc]{$_{\times 1}$}}
\put(104.33,100.33){\line(0,-1){6.00}}
\put(104.33,100.33){\line(1,0){14.00}}
\put(104.33,94.33){\line(1,0){14.00}}
\put(111.33,86.05){\makebox(0,0)[cc]{$_{162345}$}}
\put(104.33,90.39){\line(0,-1){6.00}}
\put(104.33,90.39){\line(1,0){14.00}}
\put(104.33,84.39){\line(1,0){14.00}}
\put(114.33,92.33){\makebox(0,0)[cc]{$_{\div 5}$}}
\put(92.33,78.59){\makebox(0,0)[cc]{$_{15234}$}}
\put(86.33,82.92){\line(0,-1){6.00}}
\put(98.33,82.92){\line(0,-1){6.00}}
\put(86.33,82.92){\line(1,0){12.00}}
\put(86.33,76.92){\line(1,0){12.00}}
\put(111.33,78.59){\makebox(0,0)[cc]{$_{652341}$}}
\put(101.33,81.48){\makebox(0,0)[cc]{$_{\times 1}$}}
\put(104.33,82.92){\line(0,-1){6.00}}
\put(104.33,82.92){\line(1,0){14.00}}
\put(104.33,76.92){\line(1,0){14.00}}
\put(95.33,92.33){\makebox(0,0)[cc]{$_{\div 4}$}}
\put(75.33,71.42){\makebox(0,0)[cc]{$_{1423}$}}
\put(70.33,75.76){\line(0,-1){6.00}}
\put(80.33,75.76){\line(0,-1){6.00}}
\put(70.33,69.76){\line(1,0){10.00}}
\put(70.33,75.76){\line(1,0){10.00}}
\put(92.33,71.42){\makebox(0,0)[cc]{$_{54231}$}}
\put(83.33,74.32){\makebox(0,0)[cc]{$_{\times 2}$}}
\put(86.33,75.76){\line(0,-1){6.00}}
\put(98.33,75.76){\line(0,-1){6.00}}
\put(86.33,75.76){\line(1,0){12.00}}
\put(86.33,69.76){\line(1,0){12.00}}
\put(111.33,71.42){\makebox(0,0)[cc]{$_{642315}$}}
\put(101.33,74.32){\makebox(0,0)[cc]{$_{\times 1}$}}
\put(104.33,75.76){\line(0,-1){6.00}}
\put(104.33,75.76){\line(1,0){14.00}}
\put(104.33,69.76){\line(1,0){14.00}}
\put(111.33,61.21){\makebox(0,0)[cc]{$_{142365}$}}
\put(104.33,65.55){\line(0,-1){6.00}}
\put(104.33,65.55){\line(1,0){14.00}}
\put(104.33,59.55){\line(1,0){14.00}}
\put(114.33,67.76){\makebox(0,0)[cc]{$_{\div 2}$}}
\put(78.33,92.33){\makebox(0,0)[cc]{$_{\div 3}$}}
\put(60.33,54.07){\makebox(0,0)[cc]{$_{132}$}}
\put(56.33,58.41){\line(1,0){8.00}}
\put(56.33,52.41){\line(1,0){8.00}}
\put(56.33,58.41){\line(0,-1){6.00}}
\put(64.33,58.41){\line(0,-1){6.00}}
\put(75.33,54.07){\makebox(0,0)[cc]{$_{4321}$}}
\put(67.33,56.97){\makebox(0,0)[cc]{$_{\times 3}$}}
\put(70.33,58.41){\line(0,-1){6.00}}
\put(80.33,58.41){\line(0,-1){6.00}}
\put(70.33,52.41){\line(1,0){10.00}}
\put(70.33,58.41){\line(1,0){10.00}}
\put(92.33,54.07){\makebox(0,0)[cc]{$_{53214}$}}
\put(83.33,56.97){\makebox(0,0)[cc]{$_{\times 2}$}}
\put(86.33,58.41){\line(0,-1){6.00}}
\put(98.33,58.41){\line(0,-1){6.00}}
\put(86.33,58.41){\line(1,0){12.00}}
\put(86.33,52.41){\line(1,0){12.00}}
\put(111.33,54.07){\makebox(0,0)[cc]{$_{632145}$}}
\put(104.33,58.41){\line(0,-1){6.00}}
\put(104.33,58.41){\line(1,0){14.00}}
\put(104.33,52.41){\line(1,0){14.00}}
\put(111.33,43.92){\makebox(0,0)[cc]{$_{132645}$}}
\put(104.33,48.26){\line(0,-1){6.00}}
\put(104.33,48.26){\line(1,0){14.00}}
\put(104.33,42.26){\line(1,0){14.00}}
\put(114.33,50.41){\makebox(0,0)[cc]{$_{\div 3}$}}
\put(92.33,36.77){\makebox(0,0)[cc]{$_{13254}$}}
\put(86.33,41.11){\line(0,-1){6.00}}
\put(98.33,41.11){\line(0,-1){6.00}}
\put(86.33,41.11){\line(1,0){12.00}}
\put(86.33,35.11){\line(1,0){12.00}}
\put(111.33,36.77){\makebox(0,0)[cc]{$_{632541}$}}
\put(101.33,39.67){\makebox(0,0)[cc]{$_{\times 1}$}}
\put(104.33,41.11){\line(0,-1){6.00}}
\put(104.33,41.11){\line(1,0){14.00}}
\put(104.33,35.11){\line(1,0){14.00}}
\put(95.33,50.41){\makebox(0,0)[cc]{$_{\div 4}$}}
\put(63.00,92.33){\makebox(0,0)[cc]{$_{\div 2}$}}
\put(35.44,98.33){\makebox(0,0)[cc]{$_{0,1}$}}
\put(47.44,98.33){\makebox(0,0)[cc]{$_{1,5}$}}
\put(60.44,98.33){\makebox(0,0)[cc]{$_{2,20}$}}
\put(75.44,98.33){\makebox(0,0)[cc]{$_{3,60}$}}
\put(92.44,98.33){\makebox(0,0)[cc]{$_{4,120}$}}
\put(111.44,98.33){\makebox(0,0)[cc]{$_{5,120}$}}
\put(111.44,88.39){\makebox(0,0)[cc]{$_{6,24}$}}
\put(92.44,80.92){\makebox(0,0)[cc]{$_{5,30}$}}
\put(111.44,80.92){\makebox(0,0)[cc]{$_{6,30}$}}
\put(75.44,73.76){\makebox(0,0)[cc]{$_{4,20}$}}
\put(92.44,73.76){\makebox(0,0)[cc]{$_{5,40}$}}
\put(111.44,73.76){\makebox(0,0)[cc]{$_{6,40}$}}
\put(111.44,63.55){\makebox(0,0)[cc]{$_{7,20}$}}
\put(60.44,56.41){\makebox(0,0)[cc]{$_{3,10}$}}
\put(75.44,56.41){\makebox(0,0)[cc]{$_{4,30}$}}
\put(92.44,56.41){\makebox(0,0)[cc]{$_{5,60}$}}
\put(111.44,56.41){\makebox(0,0)[cc]{$_{7,60}$}}
\put(111.44,46.26){\makebox(0,0)[cc]{$_{7,20}$}}
\put(92.44,39.11){\makebox(0,0)[cc]{$_{6,15}$}}
\put(111.44,39.11){\makebox(0,0)[cc]{$_{7,15}$}}
\put(110.76,66.26){\line(-1,0){7.06}}
\put(118.19,100.33){\line(0,-1){6.00}}
\put(118.19,90.39){\line(0,-1){6.00}}
\put(118.19,82.92){\line(0,-1){6.00}}
\put(118.19,75.76){\line(0,-1){6.00}}
\put(118.19,65.55){\line(0,-1){6.00}}
\put(118.19,58.41){\line(0,-1){6.00}}
\put(118.19,48.26){\line(0,-1){6.00}}
\put(118.19,41.11){\line(0,-1){6.00}}
\put(121.06,44.82){\makebox(0,0)[cc]{$_A$}}
\put(121.06,62.15){\makebox(0,0)[cc]{$_A$}}
\put(103.66,66.33){\line(0,-1){7.33}}
\put(103.66,59.00){\line(1,0){15.33}}
\put(119.00,59.00){\line(0,1){7.33}}
\put(119.00,66.33){\line(-1,0){2.67}}
\put(110.76,48.93){\line(-1,0){7.06}}
\put(103.66,49.00){\line(0,-1){7.33}}
\put(103.66,41.67){\line(1,0){15.33}}
\put(119.00,41.67){\line(0,1){7.33}}
\put(119.00,49.00){\line(-1,0){2.67}}
\put(18.00,76.33){\makebox(0,0)[cc]{\Large $\Lambda_3$}}
\put(17.00,55.67){\makebox(0,0)[cc]{\Large $\Lambda_4$}}
\put(23.00,29.34){\makebox(0,0)[cc]{\Large $\Lambda_5$}}
\put(55.33,89.34){\makebox(0,0)[cc]{\Large $\Lambda_6$}}
\put(6.00,84.33){\vector(1,0){6.00}}
\put(18.00,84.33){\vector(1,0){6.00}}
\put(28.00,81.33){\vector(0,-1){4.00}}
\put(6.00,65.00){\vector(1,0){6.00}}
\put(18.00,65.00){\vector(1,0){6.00}}
\put(32.00,65.00){\vector(1,0){6.00}}
\put(43.00,62.00){\vector(0,-1){4.00}}
\put(28.00,62.00){\vector(0,-1){11.33}}
\put(32.00,47.67){\vector(1,0){6.00}}
\put(6.00,37.33){\vector(1,0){6.00}}
\put(18.00,37.33){\vector(1,0){6.00}}
\put(32.00,37.33){\vector(1,0){6.00}}
\put(48.00,37.33){\vector(1,0){6.00}}
\put(60.00,34.33){\vector(0,-1){4.00}}
\put(43.00,34.33){\vector(0,-1){11.33}}
\put(48.00,20.00){\vector(1,0){6.00}}
\put(28.00,34.33){\vector(0,-1){18.67}}
\put(32.00,12.67){\vector(1,0){6.00}}
\put(48.00,12.67){\vector(1,0){6.00}}
\put(60.00,9.67){\vector(0,-1){4.00}}
\put(60.33,94.33){\vector(0,-1){36.00}}
\put(64.33,55.33){\vector(1,0){6.00}}
\put(80.33,55.33){\vector(1,0){5.67}}
\put(98.33,55.33){\vector(1,0){6.00}}
\put(111.67,52.33){\vector(0,-1){4.00}}
\put(92.33,52.33){\vector(0,-1){11.33}}
\put(98.33,38.00){\vector(1,0){6.00}}
\put(75.33,94.33){\vector(0,-1){18.67}}
\put(80.33,72.67){\vector(1,0){6.00}}
\put(98.33,72.67){\vector(1,0){6.00}}
\put(111.67,69.67){\vector(0,-1){4.00}}
\put(92.33,94.33){\vector(0,-1){11.33}}
\put(98.33,80.00){\vector(1,0){6.00}}
\put(38.33,97.33){\vector(1,0){6.00}}
\put(50.33,97.33){\vector(1,0){6.00}}
\put(64.33,97.33){\vector(1,0){6.00}}
\put(80.33,97.33){\vector(1,0){6.00}}
\put(98.33,97.33){\vector(1,0){6.00}}
\put(111.33,94.33){\vector(0,-1){4.00}}
\put(101.33,57.00){\makebox(0,0)[cc]{$_{\times 1}$}}
\put(3.00,95.70){\makebox(0,0)[cc]{$_{1}$}}
\put(15.00,95.70){\makebox(0,0)[cc]{$_{21}$}}
\put(9.00,98.59){\makebox(0,0)[cc]{$_{\times 1}$}}
\put(0.00,100.03){\line(1,0){6.00}}
\put(12.00,100.04){\line(1,0){6.00}}
\put(0.00,94.03){\line(1,0){6.00}}
\put(12.00,94.04){\line(1,0){6.00}}
\put(0.00,100.03){\line(0,-1){6.00}}
\put(6.00,100.04){\line(0,-1){6.00}}
\put(12.00,100.04){\line(0,-1){6.00}}
\put(18.00,100.04){\line(0,-1){6.00}}
\put(3.11,98.04){\makebox(0,0)[cc]{$_{0,1}$}}
\put(15.11,98.04){\makebox(0,0)[cc]{$_{1,2}$}}
\put(9.00,92.00){\makebox(0,0)[cc]{\Large $\Lambda_2$}}
\put(6.00,97.00){\vector(1,0){6.00}}
\end{picture}
\caption{Representations of $\Lambda_n$, for $n=2,3,4,5,6$.}
\end{figure}

The {\it string length} $\lambda(\Sigma(u))$ of $u$ is defined as the number of
entries ($\leq n$) of $\Sigma(u)$.
Given an arc $e$ of $\Lambda_n$, let $u_e$ and $u^e$ be the tail and the
head of $e$, respectively.  The two types of arcs in $\Lambda_n$ are selected as follows:

{\bf (1)} arcs $e$ with
$\lambda(\Sigma(u^e))=1+\lambda(\Sigma(u_e))$, as shown in Figures 1--3,
indicated with a multiplicative operator $\times m_e$, where
$c(u^e)=c(u_e)\times m_e$, noticing that $\sigma_1(u^e)\neq 1$;

\newpage

{\bf (2)} the remaining arcs $f$,
indicated with a divisive operator $\div d_f$ determined by $c(u^f)=c(u_f)\div d_f$, noticing that $\sigma_1(u^f)=1$ and that there is not an arc $e$ of type (1) with $u_e=u^f$ and $u^e=u_f$.

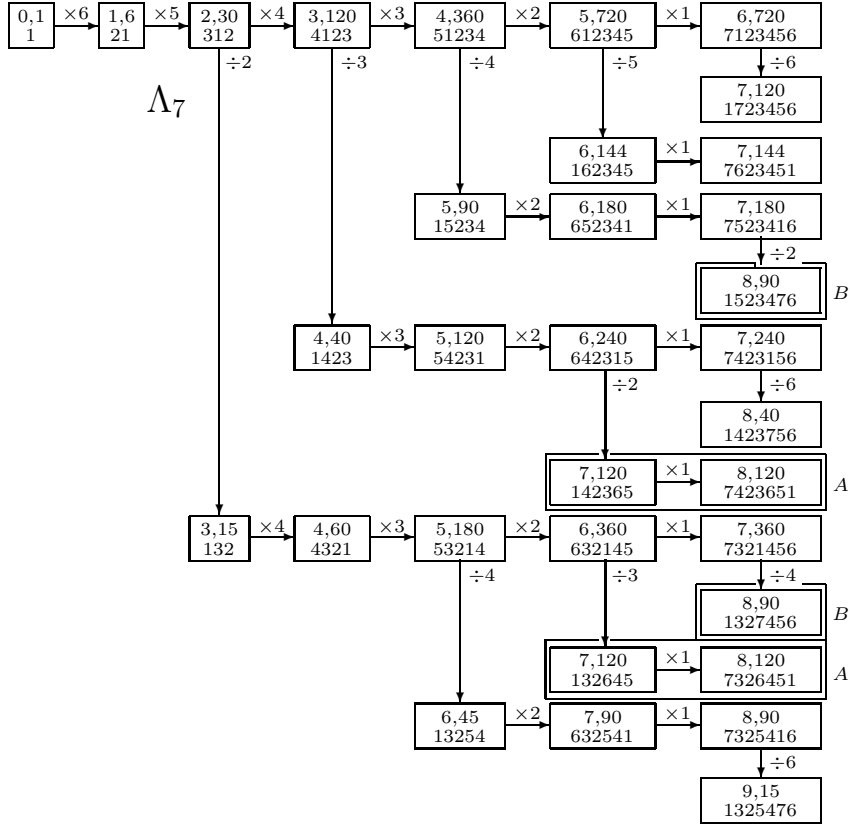
\begin{figure}[htp]
\unitlength=1.00mm
\special{em:linewidth 0.4pt}
\linethickness{0.4pt}
\begin{picture}(119.73,109.09)
\put(12.00,104.75){\makebox(0,0)[cc]{$_{1}$}}
\put(24.00,104.75){\makebox(0,0)[cc]{$_{21}$}}
\put(18.00,107.64){\makebox(0,0)[cc]{$_{\times 6}$}}
\put(37.00,104.75){\makebox(0,0)[cc]{$_{312}$}}
\put(30.00,107.64){\makebox(0,0)[cc]{$_{\times 5}$}}
\put(9.00,109.09){\line(1,0){6.00}}
\put(21.00,109.09){\line(1,0){6.00}}
\put(33.00,109.09){\line(1,0){8.00}}
\put(9.00,103.09){\line(1,0){6.00}}
\put(21.00,103.09){\line(1,0){6.00}}
\put(33.00,103.09){\line(1,0){8.00}}
\put(9.00,109.09){\line(0,-1){6.00}}
\put(15.00,109.09){\line(0,-1){6.00}}
\put(21.00,109.09){\line(0,-1){6.00}}
\put(27.00,109.09){\line(0,-1){6.00}}
\put(33.00,109.09){\line(0,-1){6.00}}
\put(41.00,109.09){\line(0,-1){6.00}}
\put(52.00,104.75){\makebox(0,0)[cc]{$_{4123}$}}
\put(44.00,107.64){\makebox(0,0)[cc]{$_{\times 4}$}}
\put(47.00,109.09){\line(0,-1){6.00}}
\put(57.00,109.09){\line(0,-1){6.00}}
\put(47.00,103.09){\line(1,0){10.00}}
\put(47.00,109.09){\line(1,0){10.00}}
\put(69.00,104.75){\makebox(0,0)[cc]{$_{51234}$}}
\put(60.00,107.64){\makebox(0,0)[cc]{$_{\times 3}$}}
\put(63.00,109.09){\line(0,-1){6.00}}
\put(75.00,109.09){\line(0,-1){6.00}}
\put(63.00,109.09){\line(1,0){12.00}}
\put(63.00,103.09){\line(1,0){12.00}}
\put(88.00,104.75){\makebox(0,0)[cc]{$_{612345}$}}
\put(78.00,107.64){\makebox(0,0)[cc]{$_{\times 2}$}}
\put(81.00,109.09){\line(0,-1){6.00}}
\put(95.00,109.09){\line(0,-1){6.00}}
\put(81.00,109.09){\line(1,0){14.00}}
\put(81.00,103.09){\line(1,0){14.00}}
\put(109.00,104.75){\makebox(0,0)[cc]{$_{7123456}$}}
\put(98.00,107.64){\makebox(0,0)[cc]{$_{\times 1}$}}
\put(101.00,109.09){\line(0,-1){6.00}}
\put(117.00,109.09){\line(0,-1){6.00}}
\put(101.00,109.09){\line(1,0){16.00}}
\put(101.00,103.09){\line(1,0){16.00}}
\put(109.00,94.93){\makebox(0,0)[cc]{$_{1723456}$}}
\put(101.00,99.26){\line(0,-1){6.00}}
\put(117.00,99.26){\line(0,-1){6.00}}
\put(101.00,99.26){\line(1,0){16.00}}
\put(101.00,93.26){\line(1,0){16.00}}
\put(112.00,101.09){\makebox(0,0)[cc]{$_{\div 6}$}}
\put(88.00,86.81){\makebox(0,0)[cc]{$_{162345}$}}
\put(81.00,91.15){\line(0,-1){6.00}}
\put(95.00,91.15){\line(0,-1){6.00}}
\put(81.00,91.15){\line(1,0){14.00}}
\put(81.00,85.15){\line(1,0){14.00}}
\put(109.00,86.81){\makebox(0,0)[cc]{$_{7623451}$}}
\put(98.00,89.71){\makebox(0,0)[cc]{$_{\times 1}$}}
\put(101.00,91.15){\line(0,-1){6.00}}
\put(117.00,91.15){\line(0,-1){6.00}}
\put(101.00,91.15){\line(1,0){16.00}}
\put(101.00,85.15){\line(1,0){16.00}}
\put(91.00,101.09){\makebox(0,0)[cc]{$_{\div 5}$}}
\put(69.00,79.35){\makebox(0,0)[cc]{$_{15234}$}}
\put(63.00,83.68){\line(0,-1){6.00}}
\put(75.00,83.68){\line(0,-1){6.00}}
\put(63.00,83.68){\line(1,0){12.00}}
\put(63.00,77.68){\line(1,0){12.00}}
\put(88.00,79.35){\makebox(0,0)[cc]{$_{652341}$}}
\put(78.00,82.24){\makebox(0,0)[cc]{$_{\times 2}$}}
\put(81.00,83.68){\line(0,-1){6.00}}
\put(95.00,83.68){\line(0,-1){6.00}}
\put(81.00,83.68){\line(1,0){14.00}}
\put(81.00,77.68){\line(1,0){14.00}}
\put(109.00,79.35){\makebox(0,0)[cc]{$_{7523416}$}}
\put(98.00,82.24){\makebox(0,0)[cc]{$_{\times 1}$}}
\put(101.00,83.68){\line(0,-1){6.00}}
\put(117.00,83.68){\line(0,-1){6.00}}
\put(101.00,83.68){\line(1,0){16.00}}
\put(101.00,77.68){\line(1,0){16.00}}
\put(109.00,69.52){\makebox(0,0)[cc]{$_{1523476}$}}
\put(101.00,73.86){\line(0,-1){6.00}}
\put(101.00,73.86){\line(1,0){16.00}}
\put(101.00,67.86){\line(1,0){16.00}}
\put(112.00,75.68){\makebox(0,0)[cc]{$_{\div 2}$}}
\put(72.00,101.09){\makebox(0,0)[cc]{$_{\div 4}$}}
\put(52.00,61.85){\makebox(0,0)[cc]{$_{1423}$}}
\put(47.00,66.19){\line(0,-1){6.00}}
\put(57.00,66.19){\line(0,-1){6.00}}
\put(47.00,60.19){\line(1,0){10.00}}
\put(47.00,66.19){\line(1,0){10.00}}
\put(69.00,61.85){\makebox(0,0)[cc]{$_{54231}$}}
\put(60.00,64.75){\makebox(0,0)[cc]{$_{\times 3}$}}
\put(63.00,66.19){\line(0,-1){6.00}}
\put(75.00,66.19){\line(0,-1){6.00}}
\put(63.00,66.19){\line(1,0){12.00}}
\put(63.00,60.19){\line(1,0){12.00}}
\put(88.00,61.85){\makebox(0,0)[cc]{$_{642315}$}}
\put(78.00,64.75){\makebox(0,0)[cc]{$_{\times 2}$}}
\put(81.00,66.19){\line(0,-1){6.00}}
\put(95.00,66.19){\line(0,-1){6.00}}
\put(81.00,66.19){\line(1,0){14.00}}
\put(81.00,60.19){\line(1,0){14.00}}
\put(109.00,61.85){\makebox(0,0)[cc]{$_{7423156}$}}
\put(98.00,64.75){\makebox(0,0)[cc]{$_{\times 1}$}}
\put(101.00,66.19){\line(0,-1){6.00}}
\put(117.00,66.19){\line(0,-1){6.00}}
\put(101.00,66.19){\line(1,0){16.00}}
\put(101.00,60.19){\line(1,0){16.00}}
\put(109.00,51.64){\makebox(0,0)[cc]{$_{1423756}$}}
\put(101.00,55.97){\line(0,-1){6.00}}
\put(117.00,55.97){\line(0,-1){6.00}}
\put(101.00,55.97){\line(1,0){16.00}}
\put(101.00,49.97){\line(1,0){16.00}}
\put(112.00,58.19){\makebox(0,0)[cc]{$_{\div 6}$}}
\put(88.00,43.97){\makebox(0,0)[cc]{$_{142365}$}}
\put(81.00,48.31){\line(0,-1){6.00}}
\put(95.00,48.31){\line(0,-1){6.00}}
\put(81.00,48.31){\line(1,0){14.00}}
\put(81.00,42.31){\line(1,0){14.00}}
\put(109.00,43.97){\makebox(0,0)[cc]{$_{7423651}$}}
\put(98.00,46.86){\makebox(0,0)[cc]{$_{\times 1}$}}
\put(101.00,48.31){\line(0,-1){6.00}}
\put(117.00,48.31){\line(0,-1){6.00}}
\put(101.00,48.31){\line(1,0){16.00}}
\put(101.00,42.31){\line(1,0){16.00}}
\put(91.00,58.19){\makebox(0,0)[cc]{$_{\div 2}$}}
\put(55.00,101.09){\makebox(0,0)[cc]{$_{\div 3}$}}
\put(37.00,36.50){\makebox(0,0)[cc]{$_{132}$}}
\put(33.00,40.84){\line(1,0){8.00}}
\put(33.00,34.84){\line(1,0){8.00}}
\put(33.00,40.84){\line(0,-1){6.00}}
\put(41.00,40.84){\line(0,-1){6.00}}
\put(52.00,36.50){\makebox(0,0)[cc]{$_{4321}$}}
\put(44.00,39.40){\makebox(0,0)[cc]{$_{\times 4}$}}
\put(47.00,40.84){\line(0,-1){6.00}}
\put(57.00,40.84){\line(0,-1){6.00}}
\put(47.00,34.84){\line(1,0){10.00}}
\put(47.00,40.84){\line(1,0){10.00}}
\put(69.00,36.50){\makebox(0,0)[cc]{$_{53214}$}}
\put(60.00,39.40){\makebox(0,0)[cc]{$_{\times 3}$}}
\put(63.00,40.84){\line(0,-1){6.00}}
\put(75.00,40.84){\line(0,-1){6.00}}
\put(63.00,40.84){\line(1,0){12.00}}
\put(63.00,34.84){\line(1,0){12.00}}
\put(88.00,36.50){\makebox(0,0)[cc]{$_{632145}$}}
\put(78.00,39.40){\makebox(0,0)[cc]{$_{\times 2}$}}
\put(81.00,40.84){\line(0,-1){6.00}}
\put(95.00,40.84){\line(0,-1){6.00}}
\put(81.00,40.84){\line(1,0){14.00}}
\put(81.00,34.84){\line(1,0){14.00}}
\put(109.00,36.50){\makebox(0,0)[cc]{$_{7321456}$}}
\put(98.00,39.40){\makebox(0,0)[cc]{$_{\times 1}$}}
\put(101.00,40.84){\line(0,-1){6.00}}
\put(117.00,40.84){\line(0,-1){6.00}}
\put(101.00,40.84){\line(1,0){16.00}}
\put(101.00,34.84){\line(1,0){16.00}}
\put(109.00,26.68){\makebox(0,0)[cc]{$_{1327456}$}}
\put(101.00,31.01){\line(0,-1){6.00}}
\put(117.00,31.01){\line(0,-1){6.00}}
\put(101.00,31.01){\line(1,0){16.00}}
\put(101.00,25.01){\line(1,0){16.00}}
\put(112.00,32.84){\makebox(0,0)[cc]{$_{\div 4}$}}
\put(88.00,19.01){\makebox(0,0)[cc]{$_{132645}$}}
\put(81.00,23.35){\line(0,-1){6.00}}
\put(95.00,23.35){\line(0,-1){6.00}}
\put(81.00,23.35){\line(1,0){14.00}}
\put(81.00,17.35){\line(1,0){14.00}}
\put(109.00,19.01){\makebox(0,0)[cc]{$_{7326451}$}}
\put(98.00,21.90){\makebox(0,0)[cc]{$_{\times 1}$}}
\put(101.00,23.35){\line(0,-1){6.00}}
\put(117.00,23.35){\line(0,-1){6.00}}
\put(101.00,23.35){\line(1,0){16.00}}
\put(101.00,17.35){\line(1,0){16.00}}
\put(91.00,32.84){\makebox(0,0)[cc]{$_{\div 3}$}}
\put(69.00,11.53){\makebox(0,0)[cc]{$_{13254}$}}
\put(63.00,15.87){\line(0,-1){6.00}}
\put(75.00,15.87){\line(0,-1){6.00}}
\put(63.00,15.87){\line(1,0){12.00}}
\put(63.00,9.87){\line(1,0){12.00}}
\put(88.00,11.53){\makebox(0,0)[cc]{$_{632541}$}}
\put(78.00,14.43){\makebox(0,0)[cc]{$_{\times 2}$}}
\put(81.00,15.87){\line(0,-1){6.00}}
\put(95.00,15.87){\line(0,-1){6.00}}
\put(81.00,15.87){\line(1,0){14.00}}
\put(81.00,9.87){\line(1,0){14.00}}
\put(109.00,11.53){\makebox(0,0)[cc]{$_{7325416}$}}
\put(98.00,14.43){\makebox(0,0)[cc]{$_{\times 1}$}}
\put(101.00,15.87){\line(0,-1){6.00}}
\put(117.00,15.87){\line(0,-1){6.00}}
\put(101.00,15.87){\line(1,0){16.00}}
\put(101.00,9.87){\line(1,0){16.00}}
\put(109.00,1.66){\makebox(0,0)[cc]{$_{1325476}$}}
\put(101.00,6.00){\line(0,-1){6.00}}
\put(117.00,6.00){\line(0,-1){6.00}}
\put(101.00,6.00){\line(1,0){16.00}}
\put(101.00,0.00){\line(1,0){16.00}}
\put(112.00,7.87){\makebox(0,0)[cc]{$_{\div 6}$}}
\put(72.00,32.84){\makebox(0,0)[cc]{$_{\div 4}$}}
\put(39.67,101.09){\makebox(0,0)[cc]{$_{\div 2}$}}
\put(12.11,107.09){\makebox(0,0)[cc]{$_{0,1}$}}
\put(24.11,107.09){\makebox(0,0)[cc]{$_{1,6}$}}
\put(37.11,107.09){\makebox(0,0)[cc]{$_{2,30}$}}
\put(52.11,107.09){\makebox(0,0)[cc]{$_{3,120}$}}
\put(69.11,107.09){\makebox(0,0)[cc]{$_{4,360}$}}
\put(88.11,107.09){\makebox(0,0)[cc]{$_{5,720}$}}
\put(109.11,107.09){\makebox(0,0)[cc]{$_{6,720}$}}
\put(109.11,97.26){\makebox(0,0)[cc]{$_{7,120}$}}
\put(88.11,89.15){\makebox(0,0)[cc]{$_{6,144}$}}
\put(109.11,89.15){\makebox(0,0)[cc]{$_{7,144}$}}
\put(69.11,81.68){\makebox(0,0)[cc]{$_{5,90}$}}
\put(88.11,81.68){\makebox(0,0)[cc]{$_{6,180}$}}
\put(109.11,81.68){\makebox(0,0)[cc]{$_{7,180}$}}
\put(109.11,71.86){\makebox(0,0)[cc]{$_{8,90}$}}
\put(52.11,64.19){\makebox(0,0)[cc]{$_{4,40}$}}
\put(69.11,64.19){\makebox(0,0)[cc]{$_{5,120}$}}
\put(88.11,64.19){\makebox(0,0)[cc]{$_{6,240}$}}
\put(109.11,64.19){\makebox(0,0)[cc]{$_{7,240}$}}
\put(109.11,53.97){\makebox(0,0)[cc]{$_{8,40}$}}
\put(88.11,46.31){\makebox(0,0)[cc]{$_{7,120}$}}
\put(109.11,46.31){\makebox(0,0)[cc]{$_{8,120}$}}
\put(37.11,38.84){\makebox(0,0)[cc]{$_{3,15}$}}
\put(52.11,38.84){\makebox(0,0)[cc]{$_{4,60}$}}
\put(69.11,38.84){\makebox(0,0)[cc]{$_{5,180}$}}
\put(88.11,38.84){\makebox(0,0)[cc]{$_{6,360}$}}
\put(109.11,38.84){\makebox(0,0)[cc]{$_{7,360}$}}
\put(109.11,29.01){\makebox(0,0)[cc]{$_{8,90}$}}
\put(88.11,21.35){\makebox(0,0)[cc]{$_{7,120}$}}
\put(109.11,21.35){\makebox(0,0)[cc]{$_{8,120}$}}
\put(69.11,13.87){\makebox(0,0)[cc]{$_{6,45}$}}
\put(88.11,13.87){\makebox(0,0)[cc]{$_{7,90}$}}
\put(109.11,13.87){\makebox(0,0)[cc]{$_{8,90}$}}
\put(109.11,4.00){\makebox(0,0)[cc]{$_{9,15}$}}
\put(117.63,41.57){\line(0,1){7.45}}
\put(117.63,49.02){\line(-1,0){28.63}}
\put(87.43,49.02){\line(-1,0){7.06}}
\put(80.37,49.02){\line(0,-1){7.45}}
\put(80.37,41.57){\line(1,0){37.25}}
\put(117.63,67.06){\line(0,1){7.45}}
\put(117.63,74.51){\line(-1,0){3.14}}
\put(108.22,73.73){\line(0,1){0.78}}
\put(108.22,74.51){\line(-1,0){7.84}}
\put(100.37,74.51){\line(0,-1){7.45}}
\put(100.37,67.06){\line(1,0){17.25}}
\put(117.63,24.31){\line(0,1){7.45}}
\put(117.63,31.76){\line(-1,0){3.14}}
\put(108.22,31.76){\line(-1,0){7.84}}
\put(100.37,31.76){\line(0,-1){7.45}}
\put(117.63,24.31){\line(-1,0){28.63}}
\put(87.43,24.31){\line(-1,0){7.06}}
\put(116.84,73.84){\line(0,-1){6.00}}
\put(119.73,70.50){\makebox(0,0)[cc]{$_B$}}
\put(119.73,45.01){\makebox(0,0)[cc]{$_A$}}
\put(119.73,27.76){\makebox(0,0)[cc]{$_B$}}
\put(119.73,19.91){\makebox(0,0)[cc]{$_A$}}
\put(80.36,16.54){\line(1,0){37.33}}
\put(80.36,16.54){\line(0,1){7.80}}
\put(117.69,16.54){\line(0,1){7.80}}
\put(15.00,106.00){\vector(1,0){6.00}}
\put(27.00,106.00){\vector(1,0){6.00}}
\put(41.00,106.00){\vector(1,0){6.00}}
\put(57.00,106.00){\vector(1,0){6.00}}
\put(75.00,106.00){\vector(1,0){6.00}}
\put(95.00,106.00){\vector(1,0){6.00}}
\put(109.00,103.00){\vector(0,-1){3.67}}
\put(88.00,103.00){\vector(0,-1){12.00}}
\put(95.00,88.00){\vector(1,0){6.00}}
\put(69.00,103.00){\vector(0,-1){19.33}}
\put(75.00,80.67){\vector(1,0){6.00}}
\put(95.00,80.67){\vector(1,0){6.00}}
\put(52.00,103.00){\vector(0,-1){36.67}}
\put(57.00,63.33){\vector(1,0){6.00}}
\put(75.00,63.33){\vector(1,0){6.00}}
\put(95.00,63.33){\vector(1,0){6.00}}
\put(109.00,60.33){\vector(0,-1){4.33}}
\put(88.33,60.33){\vector(0,-1){12.00}}
\put(95.00,45.33){\vector(1,0){6.00}}
\put(37.00,103.00){\vector(0,-1){62.00}}
\put(41.00,38.00){\vector(1,0){6.00}}
\put(57.00,38.00){\vector(1,0){6.00}}
\put(75.00,38.00){\vector(1,0){6.00}}
\put(95.00,38.00){\vector(1,0){6.00}}
\put(109.00,35.00){\vector(0,-1){4.00}}
\put(88.33,35.00){\vector(0,-1){11.67}}
\put(95.00,20.33){\vector(1,0){6.00}}
\put(69.00,35.00){\vector(0,-1){19.00}}
\put(75.00,13.00){\vector(1,0){6.00}}
\put(95.00,13.00){\vector(1,0){6.00}}
\put(109.00,10.00){\vector(0,-1){4.00}}
\put(109.00,78.01){\vector(0,-1){3.68}}
\put(30.00,96.00){\makebox(0,0)[cc]{\Large $\Lambda_7$}}
\end{picture}
\caption{Representation of $\Lambda_7$.}
\end{figure}

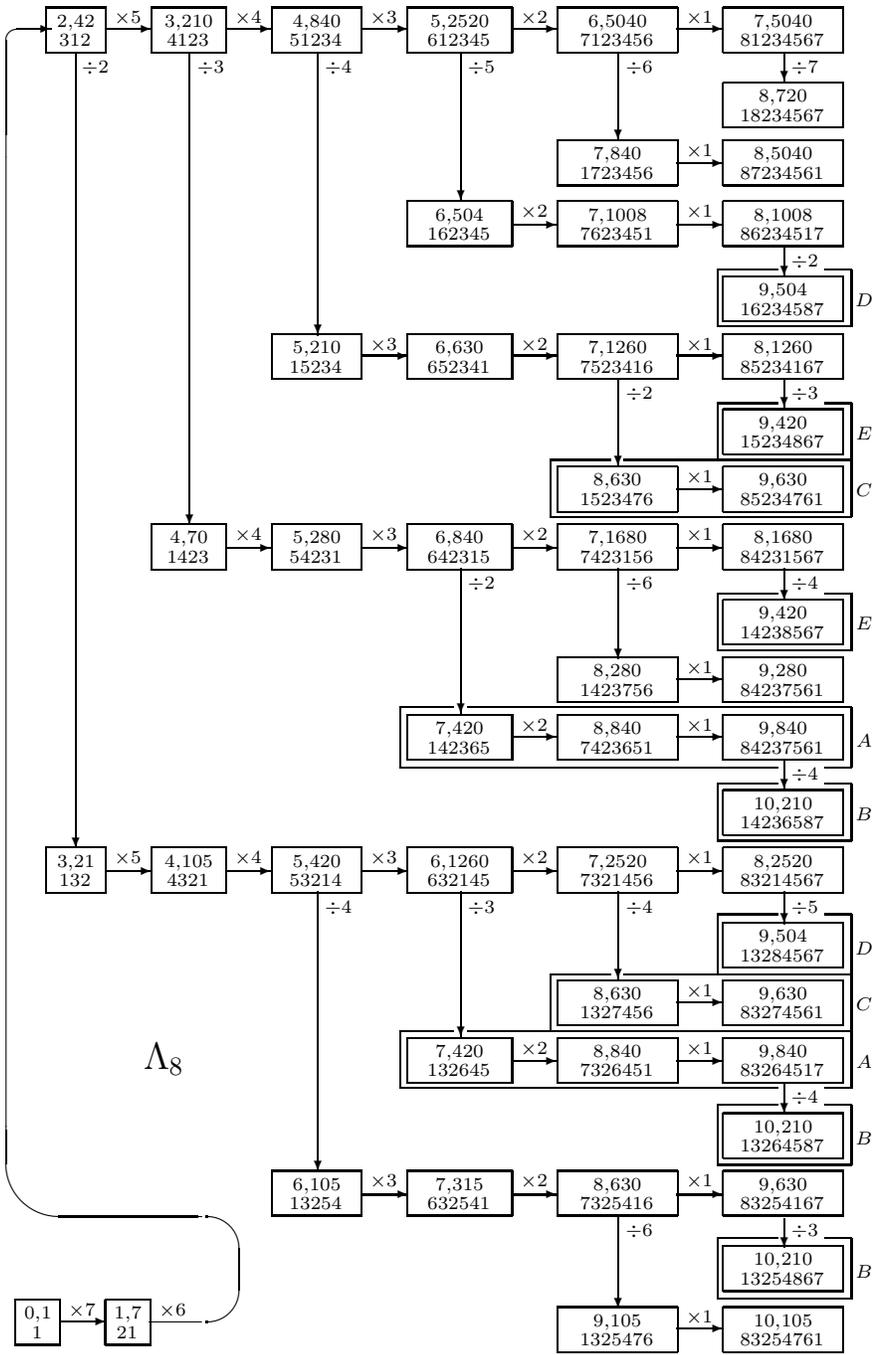
\begin{figure}[htp]
\unitlength=1.00mm
\special{em:linewidth 0.4pt}
\linethickness{0.4pt}
\begin{picture}(114.18,178.89)
\put(9.29,174.55){\makebox(0,0)[cc]{$_{312}$}}
\put(5.29,178.89){\line(1,0){8.00}}
\put(5.29,172.89){\line(1,0){8.00}}
\put(5.29,178.89){\line(0,-1){6.00}}
\put(13.29,178.89){\line(0,-1){6.00}}
\put(24.29,174.55){\makebox(0,0)[cc]{$_{4123}$}}
\put(16.29,177.44){\makebox(0,0)[cc]{$_{\times 5}$}}
\put(19.29,178.89){\line(0,-1){6.00}}
\put(29.29,178.89){\line(0,-1){6.00}}
\put(19.29,172.89){\line(1,0){10.00}}
\put(19.29,178.89){\line(1,0){10.00}}
\put(41.29,174.55){\makebox(0,0)[cc]{$_{51234}$}}
\put(32.29,177.44){\makebox(0,0)[cc]{$_{\times 4}$}}
\put(35.29,178.89){\line(0,-1){6.00}}
\put(47.29,178.89){\line(0,-1){6.00}}
\put(35.29,178.89){\line(1,0){12.00}}
\put(35.29,172.89){\line(1,0){12.00}}
\put(60.29,174.55){\makebox(0,0)[cc]{$_{612345}$}}
\put(50.29,177.44){\makebox(0,0)[cc]{$_{\times 3}$}}
\put(53.29,178.89){\line(0,-1){6.00}}
\put(67.29,178.89){\line(0,-1){6.00}}
\put(53.29,178.89){\line(1,0){14.00}}
\put(53.29,172.89){\line(1,0){14.00}}
\put(81.29,174.55){\makebox(0,0)[cc]{$_{7123456}$}}
\put(70.29,177.44){\makebox(0,0)[cc]{$_{\times 2}$}}
\put(73.29,178.89){\line(0,-1){6.00}}
\put(89.29,178.89){\line(0,-1){6.00}}
\put(73.29,178.89){\line(1,0){16.00}}
\put(73.29,172.89){\line(1,0){16.00}}
\put(103.29,174.55){\makebox(0,0)[cc]{$_{81234567}$}}
\put(92.29,177.44){\makebox(0,0)[cc]{$_{\times 1}$}}
\put(95.29,178.89){\line(0,-1){6.00}}
\put(111.29,178.89){\line(0,-1){6.00}}
\put(95.29,178.89){\line(1,0){16.00}}
\put(95.29,172.89){\line(1,0){16.00}}
\put(103.29,164.55){\makebox(0,0)[cc]{$_{18234567}$}}
\put(95.29,168.89){\line(0,-1){6.00}}
\put(111.29,168.89){\line(0,-1){6.00}}
\put(95.29,168.89){\line(1,0){16.00}}
\put(95.29,162.89){\line(1,0){16.00}}
\put(106.29,170.89){\makebox(0,0)[cc]{$_{\div 7}$}}
\put(81.29,156.89){\makebox(0,0)[cc]{$_{1723456}$}}
\put(73.29,161.22){\line(0,-1){6.00}}
\put(89.29,161.22){\line(0,-1){6.00}}
\put(73.29,161.22){\line(1,0){16.00}}
\put(73.29,155.22){\line(1,0){16.00}}
\put(103.29,156.89){\makebox(0,0)[cc]{$_{87234561}$}}
\put(92.29,159.78){\makebox(0,0)[cc]{$_{\times 1}$}}
\put(95.29,161.22){\line(0,-1){6.00}}
\put(111.29,161.22){\line(0,-1){6.00}}
\put(95.29,161.22){\line(1,0){16.00}}
\put(95.29,155.22){\line(1,0){16.00}}
\put(84.29,170.89){\makebox(0,0)[cc]{$_{\div 6}$}}
\put(60.29,148.77){\makebox(0,0)[cc]{$_{162345}$}}
\put(53.29,153.11){\line(0,-1){6.00}}
\put(67.29,153.11){\line(0,-1){6.00}}
\put(53.29,153.11){\line(1,0){14.00}}
\put(53.29,147.11){\line(1,0){14.00}}
\put(81.29,148.77){\makebox(0,0)[cc]{$_{7623451}$}}
\put(70.29,151.67){\makebox(0,0)[cc]{$_{\times 2}$}}
\put(73.29,153.11){\line(0,-1){6.00}}
\put(89.29,153.11){\line(0,-1){6.00}}
\put(73.29,153.11){\line(1,0){16.00}}
\put(73.29,147.11){\line(1,0){16.00}}
\put(103.29,148.77){\makebox(0,0)[cc]{$_{86234517}$}}
\put(92.29,151.67){\makebox(0,0)[cc]{$_{\times 1}$}}
\put(95.29,153.11){\line(0,-1){6.00}}
\put(111.29,153.11){\line(0,-1){6.00}}
\put(95.29,153.11){\line(1,0){16.00}}
\put(95.29,147.11){\line(1,0){16.00}}
\put(103.29,138.77){\makebox(0,0)[cc]{$_{16234587}$}}
\put(95.29,143.11){\line(0,-1){6.00}}
\put(111.29,143.11){\line(0,-1){6.00}}
\put(95.29,143.11){\line(1,0){16.00}}
\put(95.29,137.11){\line(1,0){16.00}}
\put(106.29,145.11){\makebox(0,0)[cc]{$_{\div 2}$}}
\put(63.29,170.89){\makebox(0,0)[cc]{$_{\div 5}$}}
\put(41.29,131.11){\makebox(0,0)[cc]{$_{15234}$}}
\put(35.29,135.44){\line(0,-1){6.00}}
\put(47.29,135.44){\line(0,-1){6.00}}
\put(35.29,135.44){\line(1,0){12.00}}
\put(35.29,129.44){\line(1,0){12.00}}
\put(60.29,131.11){\makebox(0,0)[cc]{$_{652341}$}}
\put(50.29,134.00){\makebox(0,0)[cc]{$_{\times 3}$}}
\put(53.29,135.44){\line(0,-1){6.00}}
\put(67.29,135.44){\line(0,-1){6.00}}
\put(53.29,135.44){\line(1,0){14.00}}
\put(53.29,129.44){\line(1,0){14.00}}
\put(81.29,131.11){\makebox(0,0)[cc]{$_{7523416}$}}
\put(70.29,134.00){\makebox(0,0)[cc]{$_{\times 2}$}}
\put(73.29,135.44){\line(0,-1){6.00}}
\put(89.29,135.44){\line(0,-1){6.00}}
\put(73.29,135.44){\line(1,0){16.00}}
\put(73.29,129.44){\line(1,0){16.00}}
\put(103.29,131.11){\makebox(0,0)[cc]{$_{85234167}$}}
\put(92.29,134.00){\makebox(0,0)[cc]{$_{\times 1}$}}
\put(95.29,135.44){\line(0,-1){6.00}}
\put(111.29,135.44){\line(0,-1){6.00}}
\put(95.29,135.44){\line(1,0){16.00}}
\put(95.29,129.44){\line(1,0){16.00}}
\put(103.29,121.11){\makebox(0,0)[cc]{$_{15234867}$}}
\put(95.29,125.44){\line(0,-1){6.00}}
\put(111.29,125.44){\line(0,-1){6.00}}
\put(95.29,125.44){\line(1,0){16.00}}
\put(95.29,119.44){\line(1,0){16.00}}
\put(106.29,127.44){\makebox(0,0)[cc]{$_{\div 3}$}}
\put(81.29,113.44){\makebox(0,0)[cc]{$_{1523476}$}}
\put(73.29,117.78){\line(0,-1){6.00}}
\put(89.29,117.78){\line(0,-1){6.00}}
\put(73.29,117.78){\line(1,0){16.00}}
\put(73.29,111.78){\line(1,0){16.00}}
\put(103.29,113.44){\makebox(0,0)[cc]{$_{85234761}$}}
\put(92.29,116.33){\makebox(0,0)[cc]{$_{\times 1}$}}
\put(95.29,117.78){\line(0,-1){6.00}}
\put(111.29,117.78){\line(0,-1){6.00}}
\put(95.29,117.78){\line(1,0){16.00}}
\put(95.29,111.78){\line(1,0){16.00}}
\put(84.29,127.44){\makebox(0,0)[cc]{$_{\div 2}$}}
\put(44.29,170.89){\makebox(0,0)[cc]{$_{\div 4}$}}
\put(24.29,105.77){\makebox(0,0)[cc]{$_{1423}$}}
\put(19.29,110.11){\line(0,-1){6.00}}
\put(29.29,110.11){\line(0,-1){6.00}}
\put(19.29,104.11){\line(1,0){10.00}}
\put(19.29,110.11){\line(1,0){10.00}}
\put(41.29,105.77){\makebox(0,0)[cc]{$_{54231}$}}
\put(32.29,108.67){\makebox(0,0)[cc]{$_{\times 4}$}}
\put(35.29,110.11){\line(0,-1){6.00}}
\put(47.29,110.11){\line(0,-1){6.00}}
\put(35.29,110.11){\line(1,0){12.00}}
\put(35.29,104.11){\line(1,0){12.00}}
\put(60.29,105.77){\makebox(0,0)[cc]{$_{642315}$}}
\put(50.29,108.67){\makebox(0,0)[cc]{$_{\times 3}$}}
\put(53.29,110.11){\line(0,-1){6.00}}
\put(67.29,110.11){\line(0,-1){6.00}}
\put(53.29,110.11){\line(1,0){14.00}}
\put(53.29,104.11){\line(1,0){14.00}}
\put(81.29,105.77){\makebox(0,0)[cc]{$_{7423156}$}}
\put(70.29,108.67){\makebox(0,0)[cc]{$_{\times 2}$}}
\put(73.29,110.11){\line(0,-1){6.00}}
\put(89.29,110.11){\line(0,-1){6.00}}
\put(73.29,110.11){\line(1,0){16.00}}
\put(73.29,104.11){\line(1,0){16.00}}
\put(103.29,105.77){\makebox(0,0)[cc]{$_{84231567}$}}
\put(92.29,108.67){\makebox(0,0)[cc]{$_{\times 1}$}}
\put(95.29,110.11){\line(0,-1){6.00}}
\put(111.29,110.11){\line(0,-1){6.00}}
\put(95.29,110.11){\line(1,0){16.00}}
\put(95.29,104.11){\line(1,0){16.00}}
\put(103.29,95.77){\makebox(0,0)[cc]{$_{14238567}$}}
\put(95.29,100.11){\line(0,-1){6.00}}
\put(111.29,100.11){\line(0,-1){6.00}}
\put(95.29,100.11){\line(1,0){16.00}}
\put(95.29,94.11){\line(1,0){16.00}}
\put(106.29,102.11){\makebox(0,0)[cc]{$_{\div 4}$}}
\put(81.29,88.11){\makebox(0,0)[cc]{$_{1423756}$}}
\put(73.29,92.44){\line(0,-1){6.00}}
\put(89.29,92.44){\line(0,-1){6.00}}
\put(73.29,92.44){\line(1,0){16.00}}
\put(73.29,86.44){\line(1,0){16.00}}
\put(103.29,88.11){\makebox(0,0)[cc]{$_{84237561}$}}
\put(92.29,91.00){\makebox(0,0)[cc]{$_{\times 1}$}}
\put(95.29,92.44){\line(0,-1){6.00}}
\put(111.29,92.44){\line(0,-1){6.00}}
\put(95.29,92.44){\line(1,0){16.00}}
\put(95.29,86.44){\line(1,0){16.00}}
\put(84.29,102.11){\makebox(0,0)[cc]{$_{\div 6}$}}
\put(60.29,80.44){\makebox(0,0)[cc]{$_{142365}$}}
\put(53.29,84.78){\line(0,-1){6.00}}
\put(67.29,84.78){\line(0,-1){6.00}}
\put(53.29,84.78){\line(1,0){14.00}}
\put(53.29,78.78){\line(1,0){14.00}}
\put(81.29,80.44){\makebox(0,0)[cc]{$_{7423651}$}}
\put(70.29,83.33){\makebox(0,0)[cc]{$_{\times 2}$}}
\put(73.29,84.78){\line(0,-1){6.00}}
\put(89.29,84.78){\line(0,-1){6.00}}
\put(73.29,84.78){\line(1,0){16.00}}
\put(73.29,78.78){\line(1,0){16.00}}
\put(103.29,80.44){\makebox(0,0)[cc]{$_{84237561}$}}
\put(92.29,83.33){\makebox(0,0)[cc]{$_{\times 1}$}}
\put(95.29,84.78){\line(0,-1){6.00}}
\put(111.29,84.78){\line(0,-1){6.00}}
\put(95.29,84.78){\line(1,0){16.00}}
\put(95.29,78.78){\line(1,0){16.00}}
\put(103.29,70.44){\makebox(0,0)[cc]{$_{14236587}$}}
\put(95.29,74.78){\line(0,-1){6.00}}
\put(111.29,74.78){\line(0,-1){6.00}}
\put(95.29,74.78){\line(1,0){16.00}}
\put(95.29,68.78){\line(1,0){16.00}}
\put(106.29,76.78){\makebox(0,0)[cc]{$_{\div 4}$}}
\put(63.29,102.11){\makebox(0,0)[cc]{$_{\div 2}$}}
\put(27.29,170.89){\makebox(0,0)[cc]{$_{\div 3}$}}
\put(9.29,62.77){\makebox(0,0)[cc]{$_{132}$}}
\put(5.29,67.11){\line(1,0){8.00}}
\put(5.29,61.11){\line(1,0){8.00}}
\put(5.29,67.11){\line(0,-1){6.00}}
\put(13.29,67.11){\line(0,-1){6.00}}
\put(24.29,62.77){\makebox(0,0)[cc]{$_{4321}$}}
\put(16.29,65.67){\makebox(0,0)[cc]{$_{\times 5}$}}
\put(19.29,67.11){\line(0,-1){6.00}}
\put(29.29,67.11){\line(0,-1){6.00}}
\put(19.29,61.11){\line(1,0){10.00}}
\put(19.29,67.11){\line(1,0){10.00}}
\put(41.29,62.77){\makebox(0,0)[cc]{$_{53214}$}}
\put(32.29,65.67){\makebox(0,0)[cc]{$_{\times 4}$}}
\put(35.29,67.11){\line(0,-1){6.00}}
\put(47.29,67.11){\line(0,-1){6.00}}
\put(35.29,67.11){\line(1,0){12.00}}
\put(35.29,61.11){\line(1,0){12.00}}
\put(60.29,62.77){\makebox(0,0)[cc]{$_{632145}$}}
\put(50.29,65.67){\makebox(0,0)[cc]{$_{\times 3}$}}
\put(53.29,67.11){\line(0,-1){6.00}}
\put(67.29,67.11){\line(0,-1){6.00}}
\put(53.29,67.11){\line(1,0){14.00}}
\put(53.29,61.11){\line(1,0){14.00}}
\put(81.29,62.77){\makebox(0,0)[cc]{$_{7321456}$}}
\put(70.29,65.67){\makebox(0,0)[cc]{$_{\times 2}$}}
\put(73.29,67.11){\line(0,-1){6.00}}
\put(89.29,67.11){\line(0,-1){6.00}}
\put(73.29,67.11){\line(1,0){16.00}}
\put(73.29,61.11){\line(1,0){16.00}}
\put(103.29,62.77){\makebox(0,0)[cc]{$_{83214567}$}}
\put(92.29,65.67){\makebox(0,0)[cc]{$_{\times 1}$}}
\put(95.29,67.11){\line(0,-1){6.00}}
\put(111.29,67.11){\line(0,-1){6.00}}
\put(95.29,67.11){\line(1,0){16.00}}
\put(95.29,61.11){\line(1,0){16.00}}
\put(103.29,52.77){\makebox(0,0)[cc]{$_{13284567}$}}
\put(95.29,57.11){\line(0,-1){6.00}}
\put(111.29,57.11){\line(0,-1){6.00}}
\put(95.29,57.11){\line(1,0){16.00}}
\put(95.29,51.11){\line(1,0){16.00}}
\put(106.29,59.11){\makebox(0,0)[cc]{$_{\div 5}$}}
\put(81.29,45.11){\makebox(0,0)[cc]{$_{1327456}$}}
\put(73.29,49.44){\line(0,-1){6.00}}
\put(89.29,49.44){\line(0,-1){6.00}}
\put(73.29,49.44){\line(1,0){16.00}}
\put(73.29,43.44){\line(1,0){16.00}}
\put(103.29,45.11){\makebox(0,0)[cc]{$_{83274561}$}}
\put(92.29,48.00){\makebox(0,0)[cc]{$_{\times 1}$}}
\put(95.29,49.44){\line(0,-1){6.00}}
\put(111.29,49.44){\line(0,-1){6.00}}
\put(95.29,49.44){\line(1,0){16.00}}
\put(95.29,43.44){\line(1,0){16.00}}
\put(84.29,59.11){\makebox(0,0)[cc]{$_{\div 4}$}}
\put(60.29,37.44){\makebox(0,0)[cc]{$_{132645}$}}
\put(53.29,41.78){\line(0,-1){6.00}}
\put(67.29,41.78){\line(0,-1){6.00}}
\put(53.29,41.78){\line(1,0){14.00}}
\put(53.29,35.78){\line(1,0){14.00}}
\put(81.29,37.44){\makebox(0,0)[cc]{$_{7326451}$}}
\put(70.29,40.33){\makebox(0,0)[cc]{$_{\times 2}$}}
\put(73.29,41.78){\line(0,-1){6.00}}
\put(89.29,41.78){\line(0,-1){6.00}}
\put(73.29,41.78){\line(1,0){16.00}}
\put(73.29,35.78){\line(1,0){16.00}}
\put(103.29,37.44){\makebox(0,0)[cc]{$_{83264517}$}}
\put(92.29,40.33){\makebox(0,0)[cc]{$_{\times 1}$}}
\put(95.29,41.78){\line(0,-1){6.00}}
\put(111.29,41.78){\line(0,-1){6.00}}
\put(95.29,41.78){\line(1,0){16.00}}
\put(95.29,35.78){\line(1,0){16.00}}
\put(103.29,27.44){\makebox(0,0)[cc]{$_{13264587}$}}
\put(95.29,31.78){\line(0,-1){6.00}}
\put(111.29,31.78){\line(0,-1){6.00}}
\put(95.29,31.78){\line(1,0){16.00}}
\put(95.29,25.78){\line(1,0){16.00}}
\put(106.29,33.78){\makebox(0,0)[cc]{$_{\div 4}$}}
\put(63.29,59.11){\makebox(0,0)[cc]{$_{\div 3}$}}
\put(41.29,19.77){\makebox(0,0)[cc]{$_{13254}$}}
\put(35.29,24.11){\line(0,-1){6.00}}
\put(47.29,24.11){\line(0,-1){6.00}}
\put(35.29,24.11){\line(1,0){12.00}}
\put(35.29,18.11){\line(1,0){12.00}}
\put(60.29,19.77){\makebox(0,0)[cc]{$_{632541}$}}
\put(50.29,22.67){\makebox(0,0)[cc]{$_{\times 3}$}}
\put(53.29,24.11){\line(0,-1){6.00}}
\put(67.29,24.11){\line(0,-1){6.00}}
\put(53.29,24.11){\line(1,0){14.00}}
\put(53.29,18.11){\line(1,0){14.00}}
\put(81.29,19.77){\makebox(0,0)[cc]{$_{7325416}$}}
\put(70.29,22.67){\makebox(0,0)[cc]{$_{\times 2}$}}
\put(73.29,24.11){\line(0,-1){6.00}}
\put(89.29,24.11){\line(0,-1){6.00}}
\put(73.29,24.11){\line(1,0){16.00}}
\put(73.29,18.11){\line(1,0){16.00}}
\put(103.29,19.77){\makebox(0,0)[cc]{$_{83254167}$}}
\put(92.29,22.67){\makebox(0,0)[cc]{$_{\times 1}$}}
\put(95.29,24.11){\line(0,-1){6.00}}
\put(111.29,24.11){\line(0,-1){6.00}}
\put(95.29,24.11){\line(1,0){16.00}}
\put(95.29,18.11){\line(1,0){16.00}}
\put(103.29,9.77){\makebox(0,0)[cc]{$_{13254867}$}}
\put(95.29,14.11){\line(0,-1){6.00}}
\put(111.29,14.11){\line(0,-1){6.00}}
\put(95.29,14.11){\line(1,0){16.00}}
\put(95.29,8.11){\line(1,0){16.00}}
\put(106.29,16.11){\makebox(0,0)[cc]{$_{\div 3}$}}
\put(81.29,1.66){\makebox(0,0)[cc]{$_{1325476}$}}
\put(73.29,6.00){\line(0,-1){6.00}}
\put(89.29,6.00){\line(0,-1){6.00}}
\put(73.29,6.00){\line(1,0){16.00}}
\put(73.29,0.00){\line(1,0){16.00}}
\put(103.29,1.66){\makebox(0,0)[cc]{$_{83254761}$}}
\put(92.29,4.56){\makebox(0,0)[cc]{$_{\times 1}$}}
\put(95.29,6.00){\line(0,-1){6.00}}
\put(111.29,6.00){\line(0,-1){6.00}}
\put(95.29,6.00){\line(1,0){16.00}}
\put(95.29,0.00){\line(1,0){16.00}}
\put(84.29,16.11){\makebox(0,0)[cc]{$_{\div 6}$}}
\put(44.29,59.11){\makebox(0,0)[cc]{$_{\div 4}$}}
\put(11.96,170.89){\makebox(0,0)[cc]{$_{\div 2}$}}
\put(9.40,176.89){\makebox(0,0)[cc]{$_{2,42}$}}
\put(24.40,176.89){\makebox(0,0)[cc]{$_{3,210}$}}
\put(41.40,176.89){\makebox(0,0)[cc]{$_{4,840}$}}
\put(60.40,176.89){\makebox(0,0)[cc]{$_{5,2520}$}}
\put(81.40,176.89){\makebox(0,0)[cc]{$_{6,5040}$}}
\put(103.40,176.89){\makebox(0,0)[cc]{$_{7,5040}$}}
\put(103.40,166.89){\makebox(0,0)[cc]{$_{8,720}$}}
\put(81.40,159.22){\makebox(0,0)[cc]{$_{7,840}$}}
\put(103.40,159.22){\makebox(0,0)[cc]{$_{8,5040}$}}
\put(60.40,151.11){\makebox(0,0)[cc]{$_{6,504}$}}
\put(81.40,151.11){\makebox(0,0)[cc]{$_{7,1008}$}}
\put(103.40,151.11){\makebox(0,0)[cc]{$_{8,1008}$}}
\put(103.40,141.11){\makebox(0,0)[cc]{$_{9,504}$}}
\put(41.40,133.44){\makebox(0,0)[cc]{$_{5,210}$}}
\put(60.40,133.44){\makebox(0,0)[cc]{$_{6,630}$}}
\put(81.40,133.44){\makebox(0,0)[cc]{$_{7,1260}$}}
\put(103.40,133.44){\makebox(0,0)[cc]{$_{8,1260}$}}
\put(103.40,123.44){\makebox(0,0)[cc]{$_{9,420}$}}
\put(81.40,115.78){\makebox(0,0)[cc]{$_{8,630}$}}
\put(103.40,115.78){\makebox(0,0)[cc]{$_{9,630}$}}
\put(24.40,108.11){\makebox(0,0)[cc]{$_{4,70}$}}
\put(41.40,108.11){\makebox(0,0)[cc]{$_{5,280}$}}
\put(60.40,108.11){\makebox(0,0)[cc]{$_{6,840}$}}
\put(81.40,108.11){\makebox(0,0)[cc]{$_{7,1680}$}}
\put(103.40,108.11){\makebox(0,0)[cc]{$_{8,1680}$}}
\put(103.40,98.11){\makebox(0,0)[cc]{$_{9,420}$}}
\put(81.40,90.44){\makebox(0,0)[cc]{$_{8,280}$}}
\put(103.40,90.44){\makebox(0,0)[cc]{$_{9,280}$}}
\put(60.40,82.78){\makebox(0,0)[cc]{$_{7,420}$}}
\put(81.40,82.78){\makebox(0,0)[cc]{$_{8,840}$}}
\put(103.40,82.78){\makebox(0,0)[cc]{$_{9,840}$}}
\put(103.40,72.78){\makebox(0,0)[cc]{$_{10,210}$}}
\put(9.40,65.11){\makebox(0,0)[cc]{$_{3,21}$}}
\put(24.40,65.11){\makebox(0,0)[cc]{$_{4,105}$}}
\put(41.40,65.11){\makebox(0,0)[cc]{$_{5,420}$}}
\put(60.40,65.11){\makebox(0,0)[cc]{$_{6,1260}$}}
\put(81.40,65.11){\makebox(0,0)[cc]{$_{7,2520}$}}
\put(103.40,65.11){\makebox(0,0)[cc]{$_{8,2520}$}}
\put(103.40,55.11){\makebox(0,0)[cc]{$_{9,504}$}}
\put(81.40,47.44){\makebox(0,0)[cc]{$_{8,630}$}}
\put(103.40,47.44){\makebox(0,0)[cc]{$_{9,630}$}}
\put(60.40,39.78){\makebox(0,0)[cc]{$_{7,420}$}}
\put(81.40,39.78){\makebox(0,0)[cc]{$_{8,840}$}}
\put(103.40,39.78){\makebox(0,0)[cc]{$_{9,840}$}}
\put(103.40,29.78){\makebox(0,0)[cc]{$_{10,210}$}}
\put(41.40,22.11){\makebox(0,0)[cc]{$_{6,105}$}}
\put(60.40,22.11){\makebox(0,0)[cc]{$_{7,315}$}}
\put(81.40,22.11){\makebox(0,0)[cc]{$_{8,630}$}}
\put(103.40,22.11){\makebox(0,0)[cc]{$_{9,630}$}}
\put(103.40,12.11){\makebox(0,0)[cc]{$_{10,210}$}}
\put(81.40,4.00){\makebox(0,0)[cc]{$_{9,105}$}}
\put(103.40,4.00){\makebox(0,0)[cc]{$_{10,105}$}}
\put(94.62,136.44){\line(0,1){7.56}}
\put(94.62,144.00){\line(1,0){8.00}}
\put(108.85,144.00){\line(1,0){3.56}}
\put(112.40,144.00){\line(0,-1){7.56}}
\put(112.40,136.44){\line(-1,0){17.78}}
\put(108.85,126.22){\line(1,0){3.56}}
\put(112.40,126.22){\line(0,-1){7.56}}
\put(112.40,118.67){\line(-1,0){17.78}}
\put(94.62,118.67){\line(0,1){7.56}}
\put(94.62,126.22){\line(1,0){8.00}}
\put(94.62,118.67){\line(-1,0){12.44}}
\put(80.40,118.67){\line(-1,0){8.00}}
\put(72.40,118.67){\line(0,-1){7.56}}
\put(72.40,111.11){\line(1,0){40.00}}
\put(112.40,111.11){\line(0,1){7.56}}
\put(108.85,100.89){\line(1,0){3.56}}
\put(112.40,100.89){\line(0,-1){7.56}}
\put(112.40,93.33){\line(-1,0){17.78}}
\put(94.62,93.33){\line(0,1){7.56}}
\put(94.62,100.89){\line(1,0){8.00}}
\put(61.29,85.78){\line(1,0){51.11}}
\put(52.40,85.78){\line(1,0){7.11}}
\put(112.40,85.78){\line(0,-1){8.00}}
\put(112.40,77.78){\line(-1,0){3.56}}
\put(102.62,77.78){\line(-1,0){50.22}}
\put(52.40,77.78){\line(0,1){8.00}}
\put(108.85,75.56){\line(1,0){3.56}}
\put(112.40,75.56){\line(0,-1){7.56}}
\put(112.40,68.00){\line(-1,0){17.78}}
\put(94.62,68.00){\line(0,1){7.56}}
\put(94.62,75.56){\line(1,0){8.00}}
\put(108.85,58.22){\line(1,0){3.56}}
\put(112.40,58.22){\line(0,-1){8.00}}
\put(112.40,50.22){\line(-1,0){17.78}}
\put(94.62,50.22){\line(0,1){8.00}}
\put(94.62,58.22){\line(1,0){8.00}}
\put(94.62,50.22){\line(1,0){17.78}}
\put(112.40,50.22){\line(0,-1){7.56}}
\put(112.40,42.67){\line(-1,0){40.00}}
\put(72.40,42.67){\line(0,1){7.56}}
\put(72.40,50.22){\line(1,0){8.00}}
\put(82.18,50.22){\line(1,0){12.44}}
\put(112.40,42.67){\line(0,-1){7.56}}
\put(112.40,35.11){\line(-1,0){3.56}}
\put(102.62,35.11){\line(-1,0){50.22}}
\put(52.40,35.11){\line(0,1){7.56}}
\put(59.51,42.67){\line(-1,0){7.11}}
\put(61.29,42.67){\line(1,0){11.11}}
\put(108.85,32.89){\line(1,0){3.56}}
\put(112.40,32.89){\line(0,-1){8.00}}
\put(112.40,24.89){\line(-1,0){17.78}}
\put(94.62,24.89){\line(0,1){8.00}}
\put(94.62,32.89){\line(1,0){8.00}}
\put(108.85,15.11){\line(1,0){3.56}}
\put(112.40,15.11){\line(0,-1){8.00}}
\put(112.40,7.11){\line(-1,0){17.78}}
\put(94.62,7.11){\line(0,1){8.00}}
\put(94.62,15.11){\line(1,0){8.00}}
\put(114.18,140.00){\makebox(0,0)[cc]{$_D$}}
\put(114.18,122.22){\makebox(0,0)[cc]{$_E$}}
\put(114.18,114.67){\makebox(0,0)[cc]{$_C$}}
\put(114.18,96.89){\makebox(0,0)[cc]{$_E$}}
\put(114.18,81.33){\makebox(0,0)[cc]{$_A$}}
\put(114.18,71.56){\makebox(0,0)[cc]{$_B$}}
\put(114.18,53.78){\makebox(0,0)[cc]{$_D$}}
\put(114.18,46.22){\makebox(0,0)[cc]{$_C$}}
\put(114.18,38.67){\makebox(0,0)[cc]{$_A$}}
\put(114.18,28.44){\makebox(0,0)[cc]{$_B$}}
\put(114.18,10.67){\makebox(0,0)[cc]{$_B$}}
\put(13.16,175.97){\vector(1,0){6.20}}
\put(29.44,175.97){\vector(1,0){5.81}}
\put(47.27,175.97){\vector(1,0){6.20}}
\put(67.43,175.97){\vector(1,0){5.81}}
\put(89.13,175.97){\vector(1,0){6.20}}
\put(103.47,172.87){\vector(0,-1){3.88}}
\put(81.38,172.87){\vector(0,-1){11.63}}
\put(89.13,158.14){\vector(1,0){6.20}}
\put(60.45,172.87){\vector(0,-1){19.77}}
\put(67.43,150.00){\vector(1,0){5.81}}
\put(89.13,150.00){\vector(1,0){6.20}}
\put(41.46,172.87){\vector(0,-1){37.60}}
\put(47.27,132.56){\vector(1,0){6.20}}
\put(67.43,132.56){\vector(1,0){5.81}}
\put(89.13,132.56){\vector(1,0){6.20}}
\put(103.47,129.46){\vector(0,-1){3.88}}
\put(81.38,129.46){\vector(0,-1){11.63}}
\put(89.13,114.73){\vector(1,0){6.20}}
\put(24.40,172.87){\vector(0,-1){62.79}}
\put(29.44,106.98){\vector(1,0){5.81}}
\put(47.27,106.98){\vector(1,0){6.20}}
\put(67.43,106.98){\vector(1,0){5.81}}
\put(89.13,106.98){\vector(1,0){6.20}}
\put(103.47,104.26){\vector(0,-1){4.26}}
\put(81.38,104.26){\vector(0,-1){12.02}}
\put(89.13,89.53){\vector(1,0){6.20}}
\put(60.45,104.26){\vector(0,-1){19.38}}
\put(67.43,81.78){\vector(1,0){5.81}}
\put(89.13,81.78){\vector(1,0){6.20}}
\put(103.47,78.68){\vector(0,-1){3.88}}
\put(9.29,172.87){\vector(0,-1){105.81}}
\put(13.16,63.95){\vector(1,0){6.20}}
\put(29.44,63.95){\vector(1,0){5.81}}
\put(47.27,63.95){\vector(1,0){6.20}}
\put(67.43,63.95){\vector(1,0){5.81}}
\put(89.13,63.95){\vector(1,0){6.20}}
\put(103.47,61.24){\vector(0,-1){4.26}}
\put(81.38,61.24){\vector(0,-1){11.63}}
\put(89.13,46.51){\vector(1,0){6.20}}
\put(60.45,61.24){\vector(0,-1){19.38}}
\put(67.43,38.76){\vector(1,0){5.81}}
\put(89.13,38.76){\vector(1,0){6.20}}
\put(103.47,35.66){\vector(0,-1){3.88}}
\put(41.46,61.24){\vector(0,-1){37.21}}
\put(47.27,20.93){\vector(1,0){6.20}}
\put(67.43,20.93){\vector(1,0){5.81}}
\put(89.13,20.93){\vector(1,0){6.20}}
\put(81.38,18.22){\vector(0,-1){12.40}}
\put(89.13,3.10){\vector(1,0){6.20}}
\put(103.47,147.04){\vector(0,-1){3.88}}
\put(103.47,17.74){\vector(0,-1){3.88}}
\put(4.21,2.61){\makebox(0,0)[cc]{$_{1}$}}
\put(16.21,2.61){\makebox(0,0)[cc]{$_{21}$}}
\put(10.21,5.50){\makebox(0,0)[cc]{$_{\times 7}$}}
\put(22.21,5.50){\makebox(0,0)[cc]{$_{\times 6}$}}
\put(1.21,6.95){\line(1,0){6.00}}
\put(13.21,6.95){\line(1,0){6.00}}
\put(1.21,0.95){\line(1,0){6.00}}
\put(13.21,0.95){\line(1,0){6.00}}
\put(1.21,6.95){\line(0,-1){6.00}}
\put(7.21,6.95){\line(0,-1){6.00}}
\put(13.21,6.95){\line(0,-1){6.00}}
\put(19.21,6.95){\line(0,-1){6.00}}
\put(4.32,4.95){\makebox(0,0)[cc]{$_{0,1}$}}
\put(16.32,4.95){\makebox(0,0)[cc]{$_{1,7}$}}
\put(7.11,4.03){\vector(1,0){6.20}}
\put(21.00,39.00){\makebox(0,0)[cc]{\Large $\Lambda_8$}}
\emline{19.00}{4.00}{1}{26.00}{4.00}{2}
\put(26.50,11.00){\oval(9.00,14.00)[r]}
\put(25.50,29.50){\oval(51.00,23.00)[lb]}
\emline{0.00}{30.00}{3}{0.00}{159.00}{4}
\emline{0.00}{159.00}{5}{0.00}{168.00}{6}
\emline{25.04}{18.00}{7}{26.03}{18.00}{8}
\emline{0.00}{30.02}{9}{0.00}{29.02}{10}
\put(1.99,162.00){\oval(3.98,27.95)[lt]}
\put(1.99,175.97){\vector(1,0){3.98}}
\end{picture}
\caption{Representation of $\Lambda_8$.}
\end{figure}

An additional requirement in the definition of $\Lambda_n$ is that
it is a rooted tree; its root is denoted $u_0=u_0^n$,  with $w(u_0)=0$, $c(u_0)=1$ and $\Sigma(u_0)=1$.

Given a maximal horizontal directed path, (or mhdp), $P$ of $\Lambda_n$, the
{\it depth} of $P$ is the number of vertical arcs of $\Lambda_n$ preceding
$P$ from $u_0$.

\examples Figures 1--3 contain the representations of $\Lambda_n$ for
$n=2,\ldots$, 8 (with the root of $\Lambda_8$ in Figure 3 squeezed on the bottom left), where pairs of encased mhdp's $U_I,V_I$,
either improper, (i.e. consisting of one vertex), or proper, and indicated with a common
capital letter $I=A,B,\ldots$ on their right, have corresponding vertex sets
$\{u_j^I\},\{v_j^I\}$ representing each a complete set of permutations with a common
1-ics, and thus having a common cardinality.

In fact, to determine the weight distribution of
$ST_n$, the Pruning Algorithm of section 3 below
will leave only one of these encased mpdh's with a common capital letter $I$, provided a
denomination $u_{i_0i_1\ldots i_{(j-2)}i_{(j-1)}}$ for each vertex of $\Lambda_{n+1}$
is given via the following inductive definition of Axiom $(j)$, for
$j=0,1\ldots\lfloor n/2\rfloor$, and exemplified in Figure 4, showing the strings
$i_0i_1\ldots i_{(j-2)}i_{(j-1)}$ in those denominations, for the vertices of $\Lambda_8$ in their positions in Figure 3.
{\bf Axiom (0):}
there is an mhdp $u_0u_1\ldots u_n$ of depth $0$ in
$\Lambda_{n+1}$.
{\bf Axiom ($j$):}
for each $u_{i_0i_1\ldots i_{(j-2)}i_{(j-1)}}$ as in property
($j-1$) with $i_{j-2}+1<i_{j-1}$, there is a vertical arc
$u_{i_0i_1\ldots i_{(j-1)}}u_{i_0i_1\ldots i_{(j-1)}i_{(j-1})}$ and an mhdp
from $u_{i_0i_1\ldots i_{(j-1)}i_{(j-1)}}$ to $u_{i_0i_1\ldots i_{j-1}n}$
whose depth  is $j$ in $\Lambda_{n+1}$.

\begin{figure}[htp]
\unitlength=0.50mm
\special{em:linewidth 0.4pt}
\linethickness{0.4pt}
\begin{picture}(193.08,178.89)
\put(63.08,178.89){\line(1,0){6.00}}
\put(75.08,178.89){\line(1,0){6.00}}
\put(87.08,178.89){\line(1,0){8.00}}
\put(63.08,172.89){\line(1,0){6.00}}
\put(75.08,172.89){\line(1,0){6.00}}
\put(87.08,172.89){\line(1,0){8.00}}
\put(63.08,178.89){\line(0,-1){6.00}}
\put(69.08,178.89){\line(0,-1){6.00}}
\put(75.08,178.89){\line(0,-1){6.00}}
\put(81.08,178.89){\line(0,-1){6.00}}
\put(87.08,178.89){\line(0,-1){6.00}}
\put(95.08,178.89){\line(0,-1){6.00}}
\put(101.08,178.89){\line(0,-1){6.00}}
\put(111.08,178.89){\line(0,-1){6.00}}
\put(101.08,172.89){\line(1,0){10.00}}
\put(101.08,178.89){\line(1,0){10.00}}
\put(117.08,178.89){\line(0,-1){6.00}}
\put(129.08,178.89){\line(0,-1){6.00}}
\put(117.08,178.89){\line(1,0){12.00}}
\put(117.08,172.89){\line(1,0){12.00}}
\put(135.08,178.89){\line(0,-1){6.00}}
\put(149.08,178.89){\line(0,-1){6.00}}
\put(135.08,178.89){\line(1,0){14.00}}
\put(135.08,172.89){\line(1,0){14.00}}
\put(155.08,178.89){\line(0,-1){6.00}}
\put(171.08,178.89){\line(0,-1){6.00}}
\put(155.08,178.89){\line(1,0){16.00}}
\put(155.08,172.89){\line(1,0){16.00}}
\put(177.08,178.89){\line(0,-1){6.00}}
\put(193.08,178.89){\line(0,-1){6.00}}
\put(177.08,178.89){\line(1,0){16.00}}
\put(177.08,172.89){\line(1,0){16.00}}
\put(177.08,168.89){\line(0,-1){6.00}}
\put(193.08,168.89){\line(0,-1){6.00}}
\put(177.08,168.89){\line(1,0){16.00}}
\put(177.08,162.89){\line(1,0){16.00}}
\put(155.08,161.22){\line(0,-1){6.00}}
\put(171.08,161.22){\line(0,-1){6.00}}
\put(155.08,161.22){\line(1,0){16.00}}
\put(155.08,155.22){\line(1,0){16.00}}
\put(177.08,161.22){\line(0,-1){6.00}}
\put(193.08,161.22){\line(0,-1){6.00}}
\put(177.08,161.22){\line(1,0){16.00}}
\put(177.08,155.22){\line(1,0){16.00}}
\put(135.08,153.11){\line(0,-1){6.00}}
\put(149.08,153.11){\line(0,-1){6.00}}
\put(135.08,153.11){\line(1,0){14.00}}
\put(135.08,147.11){\line(1,0){14.00}}
\put(155.08,153.11){\line(0,-1){6.00}}
\put(171.08,153.11){\line(0,-1){6.00}}
\put(155.08,153.11){\line(1,0){16.00}}
\put(155.08,147.11){\line(1,0){16.00}}
\put(177.08,153.11){\line(0,-1){6.00}}
\put(193.08,153.11){\line(0,-1){6.00}}
\put(177.08,153.11){\line(1,0){16.00}}
\put(177.08,147.11){\line(1,0){16.00}}
\put(177.08,143.11){\line(0,-1){6.00}}
\put(193.08,143.11){\line(0,-1){6.00}}
\put(177.08,143.11){\line(1,0){16.00}}
\put(177.08,137.11){\line(1,0){16.00}}
\put(117.08,135.44){\line(0,-1){6.00}}
\put(129.08,135.44){\line(0,-1){6.00}}
\put(117.08,135.44){\line(1,0){12.00}}
\put(117.08,129.44){\line(1,0){12.00}}
\put(135.08,135.44){\line(0,-1){6.00}}
\put(149.08,135.44){\line(0,-1){6.00}}
\put(135.08,135.44){\line(1,0){14.00}}
\put(135.08,129.44){\line(1,0){14.00}}
\put(155.08,135.44){\line(0,-1){6.00}}
\put(171.08,135.44){\line(0,-1){6.00}}
\put(155.08,135.44){\line(1,0){16.00}}
\put(155.08,129.44){\line(1,0){16.00}}
\put(177.08,135.44){\line(0,-1){6.00}}
\put(193.08,135.44){\line(0,-1){6.00}}
\put(177.08,135.44){\line(1,0){16.00}}
\put(177.08,129.44){\line(1,0){16.00}}
\put(177.08,125.44){\line(0,-1){6.00}}
\put(193.08,125.44){\line(0,-1){6.00}}
\put(177.08,125.44){\line(1,0){16.00}}
\put(177.08,119.44){\line(1,0){16.00}}
\put(155.08,117.78){\line(0,-1){6.00}}
\put(171.08,117.78){\line(0,-1){6.00}}
\put(155.08,117.78){\line(1,0){16.00}}
\put(155.08,111.78){\line(1,0){16.00}}
\put(177.08,117.78){\line(0,-1){6.00}}
\put(193.08,117.78){\line(0,-1){6.00}}
\put(177.08,117.78){\line(1,0){16.00}}
\put(177.08,111.78){\line(1,0){16.00}}
\put(101.08,110.11){\line(0,-1){6.00}}
\put(111.08,110.11){\line(0,-1){6.00}}
\put(101.08,104.11){\line(1,0){10.00}}
\put(101.08,110.11){\line(1,0){10.00}}
\put(117.08,110.11){\line(0,-1){6.00}}
\put(129.08,110.11){\line(0,-1){6.00}}
\put(117.08,110.11){\line(1,0){12.00}}
\put(117.08,104.11){\line(1,0){12.00}}
\put(135.08,110.11){\line(0,-1){6.00}}
\put(149.08,110.11){\line(0,-1){6.00}}
\put(135.08,110.11){\line(1,0){14.00}}
\put(135.08,104.11){\line(1,0){14.00}}
\put(155.08,110.11){\line(0,-1){6.00}}
\put(171.08,110.11){\line(0,-1){6.00}}
\put(155.08,110.11){\line(1,0){16.00}}
\put(155.08,104.11){\line(1,0){16.00}}
\put(177.08,110.11){\line(0,-1){6.00}}
\put(193.08,110.11){\line(0,-1){6.00}}
\put(177.08,110.11){\line(1,0){16.00}}
\put(177.08,104.11){\line(1,0){16.00}}
\put(177.08,100.11){\line(0,-1){6.00}}
\put(193.08,100.11){\line(0,-1){6.00}}
\put(177.08,100.11){\line(1,0){16.00}}
\put(177.08,94.11){\line(1,0){16.00}}
\put(155.08,92.44){\line(0,-1){6.00}}
\put(171.08,92.44){\line(0,-1){6.00}}
\put(155.08,92.44){\line(1,0){16.00}}
\put(155.08,86.44){\line(1,0){16.00}}
\put(177.08,92.44){\line(0,-1){6.00}}
\put(193.08,92.44){\line(0,-1){6.00}}
\put(177.08,92.44){\line(1,0){16.00}}
\put(177.08,86.44){\line(1,0){16.00}}
\put(135.08,84.78){\line(0,-1){6.00}}
\put(149.08,84.78){\line(0,-1){6.00}}
\put(135.08,84.78){\line(1,0){14.00}}
\put(135.08,78.78){\line(1,0){14.00}}
\put(155.08,84.78){\line(0,-1){6.00}}
\put(171.08,84.78){\line(0,-1){6.00}}
\put(155.08,84.78){\line(1,0){16.00}}
\put(155.08,78.78){\line(1,0){16.00}}
\put(177.08,84.78){\line(0,-1){6.00}}
\put(193.08,84.78){\line(0,-1){6.00}}
\put(177.08,84.78){\line(1,0){16.00}}
\put(177.08,78.78){\line(1,0){16.00}}
\put(177.08,74.78){\line(0,-1){6.00}}
\put(193.08,74.78){\line(0,-1){6.00}}
\put(177.08,74.78){\line(1,0){16.00}}
\put(177.08,68.78){\line(1,0){16.00}}
\put(87.08,67.11){\line(1,0){8.00}}
\put(87.08,61.11){\line(1,0){8.00}}
\put(87.08,67.11){\line(0,-1){6.00}}
\put(95.08,67.11){\line(0,-1){6.00}}
\put(101.08,67.11){\line(0,-1){6.00}}
\put(111.08,67.11){\line(0,-1){6.00}}
\put(101.08,61.11){\line(1,0){10.00}}
\put(101.08,67.11){\line(1,0){10.00}}
\put(117.08,67.11){\line(0,-1){6.00}}
\put(129.08,67.11){\line(0,-1){6.00}}
\put(117.08,67.11){\line(1,0){12.00}}
\put(117.08,61.11){\line(1,0){12.00}}
\put(135.08,67.11){\line(0,-1){6.00}}
\put(149.08,67.11){\line(0,-1){6.00}}
\put(135.08,67.11){\line(1,0){14.00}}
\put(135.08,61.11){\line(1,0){14.00}}
\put(155.08,67.11){\line(0,-1){6.00}}
\put(171.08,67.11){\line(0,-1){6.00}}
\put(155.08,67.11){\line(1,0){16.00}}
\put(155.08,61.11){\line(1,0){16.00}}
\put(177.08,67.11){\line(0,-1){6.00}}
\put(193.08,67.11){\line(0,-1){6.00}}
\put(177.08,67.11){\line(1,0){16.00}}
\put(177.08,61.11){\line(1,0){16.00}}
\put(177.08,57.11){\line(0,-1){6.00}}
\put(193.08,57.11){\line(0,-1){6.00}}
\put(177.08,57.11){\line(1,0){16.00}}
\put(177.08,51.51){\line(1,0){16.00}}
\put(155.08,49.44){\line(0,-1){6.00}}
\put(171.08,49.44){\line(0,-1){6.00}}
\put(155.08,49.44){\line(1,0){16.00}}
\put(155.08,43.44){\line(1,0){16.00}}
\put(177.08,49.44){\line(0,-1){6.00}}
\put(193.08,49.44){\line(0,-1){6.00}}
\put(177.08,49.44){\line(1,0){16.00}}
\put(177.08,43.44){\line(1,0){16.00}}
\put(135.08,41.78){\line(0,-1){6.00}}
\put(149.08,41.78){\line(0,-1){6.00}}
\put(135.08,41.78){\line(1,0){14.00}}
\put(135.08,35.78){\line(1,0){14.00}}
\put(155.08,41.78){\line(0,-1){6.00}}
\put(171.08,41.78){\line(0,-1){6.00}}
\put(155.08,41.78){\line(1,0){16.00}}
\put(155.08,35.78){\line(1,0){16.00}}
\put(177.08,41.78){\line(0,-1){6.00}}
\put(193.08,41.78){\line(0,-1){6.00}}
\put(177.08,41.78){\line(1,0){16.00}}
\put(177.08,35.78){\line(1,0){16.00}}
\put(177.08,31.78){\line(0,-1){6.00}}
\put(193.08,31.78){\line(0,-1){6.00}}
\put(177.08,31.78){\line(1,0){16.00}}
\put(177.08,25.78){\line(1,0){16.00}}
\put(117.08,24.11){\line(0,-1){6.00}}
\put(129.08,24.11){\line(0,-1){6.00}}
\put(117.08,24.11){\line(1,0){12.00}}
\put(117.08,18.11){\line(1,0){12.00}}
\put(135.08,24.11){\line(0,-1){6.00}}
\put(149.08,24.11){\line(0,-1){6.00}}
\put(135.08,24.11){\line(1,0){14.00}}
\put(135.08,18.11){\line(1,0){14.00}}
\put(155.08,24.11){\line(0,-1){6.00}}
\put(171.08,24.11){\line(0,-1){6.00}}
\put(155.08,24.11){\line(1,0){16.00}}
\put(155.08,18.11){\line(1,0){16.00}}
\put(177.08,24.11){\line(0,-1){6.00}}
\put(193.08,24.11){\line(0,-1){6.00}}
\put(177.08,24.11){\line(1,0){16.00}}
\put(177.08,18.11){\line(1,0){16.00}}
\put(177.08,14.11){\line(0,-1){6.00}}
\put(193.08,14.11){\line(0,-1){6.00}}
\put(177.08,14.11){\line(1,0){16.00}}
\put(177.08,8.11){\line(1,0){16.00}}
\put(155.08,6.00){\line(0,-1){6.00}}
\put(171.08,6.00){\line(0,-1){6.00}}
\put(155.08,6.00){\line(1,0){16.00}}
\put(155.08,0.00){\line(1,0){16.00}}
\put(177.08,6.00){\line(0,-1){6.00}}
\put(193.08,6.00){\line(0,-1){6.00}}
\put(177.08,6.00){\line(1,0){16.00}}
\put(177.08,0.00){\line(1,0){16.00}}
\put(66.19,176.11){\makebox(0,0)[cc]{${_0}$}}
\put(78.19,176.11){\makebox(0,0)[cc]{${_1}$}}
\put(91.19,176.11){\makebox(0,0)[cc]{${_2}$}}
\put(106.19,176.11){\makebox(0,0)[cc]{${_3}$}}
\put(123.19,176.11){\makebox(0,0)[cc]{${_4}$}}
\put(142.19,176.11){\makebox(0,0)[cc]{${_5}$}}
\put(163.19,176.11){\makebox(0,0)[cc]{${_6}$}}
\put(185.19,176.11){\makebox(0,0)[cc]{${_7}$}}
\put(185.19,166.11){\makebox(0,0)[cc]{${_{77}}$}}
\put(163.19,158.44){\makebox(0,0)[cc]{${_{66}}$}}
\put(185.19,158.44){\makebox(0,0)[cc]{${_{67}}$}}
\put(142.19,150.33){\makebox(0,0)[cc]{${_{55}}$}}
\put(163.19,150.33){\makebox(0,0)[cc]{${_{56}}$}}
\put(185.19,150.33){\makebox(0,0)[cc]{${_{57}}$}}
\put(185.19,140.33){\makebox(0,0)[cc]{${_{577}}$}}
\put(123.19,132.66){\makebox(0,0)[cc]{${_{44}}$}}
\put(142.19,132.66){\makebox(0,0)[cc]{${_{45}}$}}
\put(163.19,132.66){\makebox(0,0)[cc]{${_{46}}$}}
\put(185.19,132.66){\makebox(0,0)[cc]{${_{47}}$}}
\put(185.19,122.66){\makebox(0,0)[cc]{${_{477}}$}}
\put(163.19,115.00){\makebox(0,0)[cc]{${_{466}}$}}
\put(185.19,115.00){\makebox(0,0)[cc]{${_{467}}$}}
\put(106.19,107.33){\makebox(0,0)[cc]{${_{33}}$}}
\put(123.19,107.33){\makebox(0,0)[cc]{${_{34}}$}}
\put(142.19,107.33){\makebox(0,0)[cc]{${_{35}}$}}
\put(163.19,107.33){\makebox(0,0)[cc]{${_{36}}$}}
\put(185.19,107.33){\makebox(0,0)[cc]{${_{37}}$}}
\put(185.19,97.33){\makebox(0,0)[cc]{${_{377}}$}}
\put(163.19,89.66){\makebox(0,0)[cc]{${_{366}}$}}
\put(185.19,89.66){\makebox(0,0)[cc]{${_{367}}$}}
\put(142.19,82.00){\makebox(0,0)[cc]{${_{355}}$}}
\put(163.19,82.00){\makebox(0,0)[cc]{${_{356}}$}}
\put(185.19,82.00){\makebox(0,0)[cc]{${_{357}}$}}
\put(185.19,72.00){\makebox(0,0)[cc]{${_{3577}}$}}
\put(91.19,64.33){\makebox(0,0)[cc]{${_{22}}$}}
\put(106.19,64.33){\makebox(0,0)[cc]{${_{23}}$}}
\put(123.19,64.33){\makebox(0,0)[cc]{${_{24}}$}}
\put(142.19,64.33){\makebox(0,0)[cc]{${_{25}}$}}
\put(163.19,64.33){\makebox(0,0)[cc]{${_{26}}$}}
\put(185.19,64.33){\makebox(0,0)[cc]{${_{27}}$}}
\put(185.19,54.33){\makebox(0,0)[cc]{${_{277}}$}}
\put(163.19,46.66){\makebox(0,0)[cc]{${_{266}}$}}
\put(185.19,46.66){\makebox(0,0)[cc]{${_{267}}$}}
\put(142.19,39.00){\makebox(0,0)[cc]{${_{255}}$}}
\put(163.19,39.00){\makebox(0,0)[cc]{${_{256}}$}}
\put(185.19,39.00){\makebox(0,0)[cc]{${_{257}}$}}
\put(185.19,29.00){\makebox(0,0)[cc]{${_{2577}}$}}
\put(123.19,21.33){\makebox(0,0)[cc]{${_{244}}$}}
\put(142.19,21.33){\makebox(0,0)[cc]{${_{245}}$}}
\put(163.19,21.33){\makebox(0,0)[cc]{${_{246}}$}}
\put(185.19,21.33){\makebox(0,0)[cc]{${_{247}}$}}
\put(185.19,11.33){\makebox(0,0)[cc]{${_{2477}}$}}
\put(163.19,3.22){\makebox(0,0)[cc]{${_{2466}}$}}
\put(185.19,3.22){\makebox(0,0)[cc]{${_{2467}}$}}
\put(68.98,175.97){\vector(1,0){6.20}}
\put(81.00,175.97){\vector(1,0){6.20}}
\put(94.95,175.97){\vector(1,0){6.20}}
\put(111.23,175.97){\vector(1,0){5.81}}
\put(129.06,175.97){\vector(1,0){6.20}}
\put(149.22,175.97){\vector(1,0){5.81}}
\put(170.92,175.97){\vector(1,0){6.20}}
\put(185.26,172.87){\vector(0,-1){3.88}}
\put(163.17,172.87){\vector(0,-1){11.63}}
\put(170.92,158.14){\vector(1,0){6.20}}
\put(142.24,172.87){\vector(0,-1){19.77}}
\put(149.22,150.00){\vector(1,0){5.81}}
\put(170.92,150.00){\vector(1,0){6.20}}
\put(123.25,172.87){\vector(0,-1){37.60}}
\put(129.06,132.56){\vector(1,0){6.20}}
\put(149.22,132.56){\vector(1,0){5.81}}
\put(170.92,132.56){\vector(1,0){6.20}}
\put(185.26,129.46){\vector(0,-1){3.88}}
\put(163.17,129.46){\vector(0,-1){11.63}}
\put(170.92,114.73){\vector(1,0){6.20}}
\put(106.19,172.87){\vector(0,-1){62.79}}
\put(111.23,106.98){\vector(1,0){5.81}}
\put(129.06,106.98){\vector(1,0){6.20}}
\put(149.22,106.98){\vector(1,0){5.81}}
\put(170.92,106.98){\vector(1,0){6.20}}
\put(185.26,104.26){\vector(0,-1){4.26}}
\put(163.17,104.26){\vector(0,-1){12.02}}
\put(170.92,89.53){\vector(1,0){6.20}}
\put(142.24,104.26){\vector(0,-1){19.38}}
\put(149.22,81.78){\vector(1,0){5.81}}
\put(170.92,81.78){\vector(1,0){6.20}}
\put(185.26,78.68){\vector(0,-1){3.88}}
\put(91.08,172.87){\vector(0,-1){105.81}}
\put(94.95,63.95){\vector(1,0){6.20}}
\put(111.23,63.95){\vector(1,0){5.81}}
\put(129.06,63.95){\vector(1,0){6.20}}
\put(149.22,63.95){\vector(1,0){5.81}}
\put(170.92,63.95){\vector(1,0){6.20}}
\put(185.26,61.24){\vector(0,-1){4.26}}
\put(163.17,61.24){\vector(0,-1){11.63}}
\put(170.92,46.51){\vector(1,0){6.20}}
\put(142.24,61.24){\vector(0,-1){19.38}}
\put(149.22,38.76){\vector(1,0){5.81}}
\put(170.92,38.76){\vector(1,0){6.20}}
\put(185.26,35.66){\vector(0,-1){3.88}}
\put(123.25,61.24){\vector(0,-1){37.21}}
\put(129.06,20.93){\vector(1,0){6.20}}
\put(149.22,20.93){\vector(1,0){5.81}}
\put(170.92,20.93){\vector(1,0){6.20}}
\put(163.17,18.22){\vector(0,-1){12.40}}
\put(170.92,3.10){\vector(1,0){6.20}}
\put(185.26,147.04){\vector(0,-1){3.88}}
\put(185.26,17.74){\vector(0,-1){3.88}}
\end{picture}
\caption{Index string representation of $\Lambda_8$.}
\end{figure}

\section{Redefinition and pruning of $\Lambda_n$}

If a vertex $u$ of $\Lambda_n$ is the tail of a horizontal (vertical) arc $e$ indicated $\times m_e$ ($\div d_e$) then we write $m_u=m_e$ ($d_u=d_e$). If $u$ is not the tail of an arc, then we write $m_u=0$, ($d_u=0$). The following redefinition allows to consider $\Lambda_n$ as a  subdigraph of $\Lambda_{n+1}$ for any $n$, so the limit $\Lambda_\infty$ of the nested sequence $\{\Lambda_n; n>0\}$ of binary directed trees makes sense:
replace the indication $\times m_u=\times m_e$ of the horizontal arc $e$ of $\Lambda_n$ departing from a tail $u$ of an arc in $\Lambda_n$ by $\bullet \ell_u$, where $\ell_u=n-m_u$ and $\bullet$ is the operation given by $c(u)\bullet\ell_u= c(u)\times(n-m_u)$; let $\Lambda_n^\bullet$ be the resulting indicated digraph; redefine $\Lambda_n=\Lambda_n^\bullet$.
The new indications for the horizontal arcs allow
now the containment of indicated subdigraphs. Incidentally, the cardinalities of the set of vertices of the resulting $\Lambda_\infty$ that are tails of arcs indicated $\bullet_i$ form a Fibonacci sequence according to the increasing values of $i=0,1,\ldots$. This is apparent from the number of vertices in the successive columns from left to right, in the representations of the $\Lambda_n$'s as in Figures 1--3, for increasing values of $n=2,3,\ldots$.

{\bf Pruning Algorithm.} The vertices $u_{i_0i_1\ldots i_j}$ of $\Lambda_n$
are treated first in the increasing order of their string lengths $j+1$ and then,
for each fixed string length $j+1$,
in the lexicographical order of their subindex strings $i_0i_1\ldots
i_j$, namely:
$$\begin{array}{c}
u_0,u_1,\ldots,u_{n-1},u_{2,2},u_{2,3},\ldots,u_{2,n-1},u_{3,3},\ldots,\\
u_{n-1,n-1},u_{2,4,4},\ldots,u_{2,4,6,6},\ldots\end{array}$$ Each such a vertex
$u=u_{i_0i_1\ldots i_j}$ has the following fields associated with it:
{\bf (1)} the notation $u=u_i^w=u_{i_0i_1\ldots u_j}^{w(u)}$, where $w(u)=w(\Sigma(u))$;
{\bf (2)} the notation $\Sigma(u)$ of the corresponding permutation of
$\{1,\ldots n\}$ associated to $u$;
{\bf (3)} the 1-ics $\Pi(u)$ of $\Sigma(u)$;
{\bf (4)} the number $\ell_u=n-m_u$;
{\bf (5-7)} either a blank in each of the three cases (5), (6) and (7),
if $u$ is the first or second vertex of an mhdp, or:
{\bf (5)} the notation $\Sigma[u]$ of the permutation obtained from $\Sigma(u)$
by permuting $\sigma_1$ and $\sigma_k=1$, ($k\neq 1$);
{\bf (6)} the 1-ics $\Pi[u]$ of $\Sigma[u]$;
{\bf (7)} a tuple $C(u)=s_1,\ldots,s_h$ composed by the
orders $s_j$ of the cycles composing $\Pi[u]$;
{\bf (8)} the number $d_u$ expressed as a product $b_ua_u$, where
{\bf (8a)} $a_u=i_j-i_{j-1}$ under the convention $i_{-1}=0$ and
{\bf (8b)} $b_u\neq 0$ if the value of item (7) above is not a
blank and the resulting tuple $C(u)$ was not present in previously treated vertices
of $\Lambda_n$; $b_u=0$, otherwise.

The Pruning Algorithm consists in determining these fields for the vertices of $\Lambda_n$, in the prescribed order. This allows a
partial reconstruction of $\Lambda_n$ in the form of the maximal subdigraph $\Lambda'_n$, which
accepts and copies all the vertices and arcs of $\Lambda_n$ into $\Lambda'_n$
except for one case: If $b_u=0$ and there is a vertical arc $e$ whose
tail is $u$, then $e$ is not copied from $\Lambda_n$ into $\Lambda'_n$; this is interpreted as the pruning of $e$ and descendant vertices and arcs (performed to avoid repetitions of mhdp's, as in the encased mpdh's having a
common capital-letter indication at their right in Figures 1--3).\qfd

Let $\rho_n$ be the relation defined on the vertex set of $ST_n$ by
$u\rho_nv$ if and only if $u$ and $v$ represent permutations with a common
1-ics. The following second redefinition of
$\Lambda_n$ allows to have its vertex set in bijective correspondence
with the family of equivalence classes of $ST_n$ under $\rho_n$, which in
turn allows to use $\Lambda_n$ in computing the weight distribution of
$ST_n$:
perform the Pruning Algorithm of $\Lambda_n$, whose output is a
maximal subdigraph
$\Lambda'_n$ of $\Lambda_n$ in which there are not pairs of mhdp's $v_0v_1\ldots
v_s$ and $v'_0v'_1\ldots v'_s$ of the same string length $s$ with corresponding
vertices $v_i$ and $v'_i$ having common 1-ics $\Pi(v_i)=\Pi(v'_i)$,
for $i=0,1,\ldots,s$; redefine $\Lambda_n=\Lambda'_n$. We still have $\Lambda_n$
as a subdigraph of $\Lambda_{n+1}$ for every $n$, so a $\Lambda_\infty$
persists.

\example
The algorithm yields the list ${\mathcal P}_9$
for $n=9$, (commas are deleted in subindices $i$ of $u_i^w=u(i,w)$ in item 1 and in the tuples $C(u)$ in item 7;
the $\Pi(u)$ and $\Pi[u]$ are shown to the right of their corresponding
$\Sigma(u)$ and $\Sigma[u]$):

$$\begin{array}{llllllll}
_{u(i,w)}&_{\Sigma(u)\Pi(u)}&_{\ell_u}&_{\Sigma[u]\Pi[u]}&_{C(u)}&_{b_ua_u}\\ \hline
^{u(0,0)}_{u(1,1)}&^{1}_{21(12)}&_{2}&& _{1}&^{00}_{01}\\
^{u(2,2)}_{u(3,3)}&^{312(132)}_{4123(1432)}&^3_{4} &^{132(32)}_{1423(432)}&^{2}_{3}&^{12}_{13}\\
^{u(4,4)}_{u(5,5)}&^{51234(15432)}_{612345(165432)}&^{5}_{6}&^{15234(5432)}_{162345(65432)}&^{4}_{5}&^{14}_{15}\\
^{u(6,6)}_{u(7,7)}&^{7123456(1765432)}_{81234567(18765432)}&^{7}_{8}&^{1723456(765432)}_{18234567(8765432)}&^{6}_{7}&^{16}_{17}\\
^{u(8,8)}_{u(22,3)}&^{912345678(198765432)}_{132(23)}&^{9}_{3}&^{192345678(98765432)}&^{8}&^{18}_{ 00}\\
^{u(23,4)}_{u(24,5)}&^{4321(14.32)}_{53214(154.32)}&^{4}_{5}&_{13254(54.32)}&_{22}&^{01}_{22}\\
^{u(25,6)}_{u(26,7)}&^{632145(1654.32)}_{7321456(17654.32)}&^{6}_{7}&^{132645(654.32)}_{1327456(7654.32)}&^{23}_{24}&^{13}_{14}\\
^{u(27,8)}_{u(28,9)}&^{83214567(187654.32)}_{932145678(1987654.32)}&^{8}_{9}&^{13284567(87654.32)}_{132945678(987654.32)}&^{25}_{26}&^{15}_{16}\\
^{u(33,4)}_{u(34,5)}&^{1423(432)}_{54231(15.432)}&^{4}_{5}&&&^{00}_{01}\\
^{u(35,6)}_{u(36,7)}&^{642315(165.432)}_{7423156(1765.432)}&^{6}_{7}&^{142365(65.432)}_{1423756(765.432)}&^{32}_{33}&^{02}_{23}\\
^{u(37,8)}_{u(38,9)}&^{84231567(18765.432)}_{942315678(198765.432)}&^{8}_{9}&^{14238567(8765.432)}_{142395678(98765.432)}&^{34}_{35}&^{14}_{15}\\
^{u(44,5)}_{u(45,6)}&^{15234(5432)}_{652341(16.5432)}&^{5}_{6}&&&^{00}_{01}\\
^{u(46,7)}_{u(47,8)}&^{7523416(176.5432)}_{85234167(1876.5432)}&^{7}_{8}&^{1523476(76.5432)}_{15234867(876.5432)}&^{42}_{43}&^{02}_{03}\\
^{u(48,9)}_{u(55,6)}&^{952341678(19876.5432)}_{162345(65432)}&^{9}_{6}&^{152349678(9876.5432)}&^{44}&^{24}_{00}\\
^{u(56,7)}_{u(57,8)}&^{7623451(17.65432)}_{86234517(187.65432)}&^{7}_{8}&_{16234587(87.65432)}&_{52}&^{01}_{02}\\
^{u(58,9)}_{u(66,7)}&^{962345178(1987.65432)}_{1723456(765432)}&^{9}_{7}&^{162345978(987.65432)}&^{53}&^{03}_{00}\\
^{u(67,8)}_{u(68,9)}&^{87234561(18.765432)}_{972345618(198.765432)}&^{8}_{9}&_{172345698(98.765432)}&_{62}&^{01}_{02}\\
^{u(77,8)}_{u(78,9)}&^{18234567(8765432)}_{982345671(19.8765432)}&^{8}_{9}&&&^{00}_{01}\\
^{u(88,9)}_{u(244,6)}&^{192345678(98765432)}_{13254(54.32)}&^{9}_{5}&&&^{00}_{00}\\
^{u(245,7)}_{u(246,8)}&^{632541(16.54.32)}_{7325416(176.54.32)}&^{6}_{7}&_{1325476(23.45.67)}&_{222}&^{01}_{32}\\
^{u(247,9)}_{u(248,{10})}&^{83254167(1876.54.32)}_{932541678(19876.54.32)}&^{8}_{9}&^{13254867(23.45.687)}_{132549678(23.45.6987)}&^{223}_{224}&^{13}_{14}\\
^{u(255,7)}_{u(256,8)}&^{132645(654.32)}_{7326451(17.654.32)}&^{6}_{7}&&&^{00}_{01}\\
^{u(257,9)}_{u(258,{10})}&^{83264517(187.654.32)}_{932645178(1987.654.32)}&^{8}_{9}&^{13264587(23.465.78)}_{132645978(23.465.798)}&^{232}_{233}&^{02}_{13}\\
^{u(266,8)}_{u(267,9)}&^{1327456(7654.32)}_{83274561(18.7654.32)}&^{7}_{8}&&&^{00}_{01}\\
^{u(268,{10})}_{u(277,9)}&^{932745618(198.7654.32)}_{13284567(87654.32)}&^{9}_{8}&^{132745698(23.4765.89)}&^{242}&^{02}_{00}\\
^{u(278,{10})}_{u(288,{10})}&^{932845671(19.87654.32)}_{132945678(987654.32)}&^{9}_{9}&&&^{01}_{00}\\
^{u(366,8)}_{u(367,9)}&^{1423756(765.432)}_{84237561(18.765.432)}&^{7}_{8}&&&^{00}_{01}\\
^{u(368,{10})}_{u(377,9)}&^{942375618(198.765.432)}_{14238567(8765.432)}&^{9}_{8}&^{142375698(243.576.89)}&^{332}&^{02}_{00}\\
^{u(378,{10})}_{u(388,{10})}&^{942385671(19.8765.432)}_{142395678(98765.432)}&^{9}_{9}&&&^{01}_{00}\\
^{u(488,{10})}_{u(2466,9)}&^{152349678(9876.5432)}_{1325476(76.54.32)}&^{9}_{7}&&&^{00}_{00}\\
^{u(2467,{10})}_{u(2468,{11})}&^{83254761(18.76.54.32)}_{932547618(198.76.54.32)}&^{8}_{9}&_{132547698(23.45.67.89)}&_{2222}&^{01}_{42}\\
^{u(2477,{10})}_{u(2478,{11})}&^{13254867(876.54.32)}_{932548671(19.876.54.32)}&^{8}_{9}&&&^{00}_{01}\\
^{u(2488,{11})}_{u(2588,{11})}&^{132549678(9876.54.32)}_{132645978(987.654.32)}&^{9}_{9}&&&^{00}_{00}\\
^{u(24688,12)}&^{u(132547698(98.76.54.32)}&^{9}&&&^{00}\\
\end{array}$$

This list ${\mathcal P}_n$ generalizes to the patterns expressed in the following
theorem. For $u=u_{i_0i_1\ldots i_j}$ in $\Lambda_n$, let
$\ell_u=\ell_{i_0i_1\ldots i_j}$, etc.

\begin{thm} Let $i_{-1}=0$ and let $t_k=i_k-i_{k-1}$, for $k=0,1,\ldots,j-1$.
Then:
{\bf(1)} the $1$-ics $C(u)$ in the penultimate
field of the line associated to a vertex $u=u_{i_0i_1\ldots i_j}$ in
${\mathcal P}_n$ is of the form
$t_0,t_1,\ldots,t_j$, where the order of the integers $t_k$ is irrelevant;
{\bf(2)} the vertices $u_{i_0i_1\ldots i_j}$ of $\Lambda_n$,
{\rm(}remaining after applying the Pruning Algorithm{\rm)},
have subindex strings $i_0i_1\ldots i_j$
completely determined by the following conditions:
$$\begin{array}{ll}
{\rm(a)}\mbox{ }0\leq i_0\leq n-1; &
{\rm(b)}\mbox{ if }j>0,\mbox{ then }2\leq i_0; \\
{\rm(c)}\mbox{ }t_k\leq t_{k+1}, \mbox{ for } k=0,\ldots,j-2; &
{\rm(d)}\mbox{ }i_{j-1}\leq i_j;
\end{array}$$
{\bf(3)} the weight $w(u)$ of a vertex $u=u_{i_0i_1\ldots i_j}$ of
$\Lambda_n$ is $w(u)=w(u_{i_0i_1\ldots i_j})=i_j+j;$
{\bf(4)} the number $\ell_u$ associated to a vertex
$u_{i_0i_1\ldots i_j}$
of $\Lambda_n$ is $\ell_u=\ell_{i_0i_1\ldots i_j}=i_j+1;$ thus, the
corresponding multiplicative factor $m_u$ is
$m_u=m_{i_0i_1\ldots i_j}=n-i_j-1;$
{\bf(5)} the divisive-operator number $d_u=b_u.a_u$ has $a_u=t_j$; moreover,
$b_u>0$ if and only if either $j=0$ and $i_0>1$ or $j>0$ and
$2\leq i_0\leq t_1\leq t_2 \leq \ldots\leq t_j$; furthermore,
if $b_u>0$, then $b_u=1$, unless $i_0=t_1=t_2=\ldots=t_j$,
in which case $b_u=j+1$.
\end{thm}

\section{The weight distribution of $ST_n$}

To compute the weight distribution of $ST_n$, a table ${\mathcal T}_n$ constructed from the resulting pruned version of $\Lambda_n$ and satisfying
the following additional conditions will be used:
{\bf (a)} the subindex strings $i_0i_1\ldots i_j$ of the
vertices $u_{i_0i_1\ldots i_j}$ of $\Lambda_n$ are distributed on columns
according to their weights;
{\bf (b)} each row is to contain the subindex strings of the vertices of an mhdp
$P$ of $\Lambda_n$, given from left to right according to the orientation of
$P$;
{\bf (c)} each mhdp is presented in lexicographical order in its containing row;
{\bf (d)} the rows of each complete set of common-depth mhdp's
are presented contiguously and in the decreasing order of their path lengths,
thus forming upper triangular matrices, because of item (a), above;
{\bf (e)} these upper triangular matrices are given from top to bottom
in the increasing order of their depths.

\example
${\mathcal T}_{11}$ is as follows, where $a=10$ and $b=11$, vertices
with  $j=2$ and $i_0=3$ previous to $366$ do not appear since they were
pruned, and one additional row should be added for the 15-th column, containing
solely the string $2468aa$, (which, for insufficient margin, remained excluded):

$$\begin{array}{|rrrrrrrrrrrrrrr|}\hline
_{0}&_{1}&_{2}&_{3} &_{4} &_{5} &_{6} &_{7} &_{8} &_{9}&_{10}&_{11}&_{12}&_{13}&_{14}\\\hline
^{0}&^{1}&^{2}&^{\;\;3}_{22} &^{\;\;4}_{23} &^{\;\;5}_{24} &^{\;\;6}_{25} &^{\;\;7}_{26} &^{\;\;8}_{27} &^{\;\;9}_{28}&^{\;\;a}_{29}&_{2a}&&&\\
    &    &    &     &^{33}&^{34}_{44}&^{35}_{45}&^{36}_{46}&^{37}_{47}&^{38}_{48}&^{39}_{49}&^{3a}_{4a}&&&\\
    &    &    &     &     &     &^{55}&^{56}_{66}&^{57}_{67}&^{58}_{68}&^{59}_{69}&^{5a}_{6a}&&&\\
    &    &    &     &     &     &     &     &^{77}&^{78}_{88}&^{79}_{89}&^{7a}_{8a}&&&\\
    &    &    &     &     &     &     &     &     &     &^{99}&^{9a}_{aa}&&&\\
\hline\end{array}$$

$$\begin{array}{|rrrrrrrrrrrrrrr|}\hline
_{0}&_{1}&_{2}&_{3} &_{4} &_{5} &_{6} &_{7} &_{8} &_{9}&_{10}&_{11}&_{12}&_{13}&_{14}\\\hline
    &    &    &     &     &     &^{244}&^{245}_{255}&^{246}_{256}&^{247}_{257}&^{248}_{258}&^{249}_{259}&^{24a}_{25a}&&\\
    &    &    &     &     &     &&&^{266}&^{267}_{277}&^{268}_{278}&^{269}_{279}&^{26a}_{27a}&&\\
    &    &    &     &     &     &&&      &            &^{288}      &^{289}_{299}&^{28a}_{29a}&&\\
    &    &    &     &     &     &&&_{366}&_{367}      &_{368}      &_{369}      &^{2aa}_{36a}&&\\
    &    &    &     &     &     &&&      &^{377}      &^{378}_{388}&^{379}_{389}&^{37a}_{38a}&&\\
    &    &    &     &     &     &&&      &            &            &^{399}      &^{39a}_{3aa}&&\\
    &    &    &     &     &     &&&      &            &^{488}&^{489}_{499}&^{48a}_{49a}&&\\
    &    &    &     &     &     &&&      &            &&&^{4aa}_{5aa}&&\\
&&&&&&&&&^{2466}&^{2467}_{2477}&^{2468}_{2478}&^{2469}_{2479}&^{246a}_{247a}&\\
&&&&&&&&&&&^{2488}&^{2489}_{2499}&^{248a}_{249a}&\\
&&&&&&&&&&&_{2588}&_{2589}&^{24aa}_{258a}&\\
&&&&&&&&&&&&^{2599}&^{259a}_{25aa}&\\
&&&&&&&&&&&&_{3699}&^{26aa}_{369a}&\\
&&&&&&&&&&&&_{24688}&^{\;\;36aa}_{24689}&_{2468a}\\
&&&&&&&&&&&&&^{24699}&^{2469a}_{246aa}\\
&&&&&&&&&&&&&&^{247aa}\\\hline
\end{array}$$

Each vertex $u=u_{i_0i_1\ldots i_j}=i_0i_1\ldots i_j$ of
$\Lambda_n$, reachable from $u_0=0$ by a path $P$, has associated
cardinality $c(u)=M/(A.B)$, where:
{\bf (a)} $M$, (respectively $A$), is the product of the numbers
$m_{i_0i_1\ldots i_j}=n-i_j-1$, (resp. $a_{i_0i_1\ldots i_j}=t_j=i_i-i_{j-1}$),
of all tails $i_0i_1\ldots i_j$
of horizontal, (resp. vertical), arcs in $P$; {\bf (b)} $B$ is the
product of all the numbers $b_{i_0i_1\ldots i_j}$ of tails
$i_0i_1\ldots i_j$ of vertical arcs in $P$ with $i_0=t_1=\ldots=t_j$.

A procedure to compute the path from $u_0$ to any given vertex
$u$ of $\Lambda_n$, performed by going backwards from $i_0i_1\ldots i_j$ to
$0$ by means of table ${\mathcal T}_n$, consists of the following steps:
{\bf (1)} set $u=i_0i_1\ldots i_j$; {\bf (2)} if $u$ is not the
first vertex of an mhdp, then go backwards through the vertices of the mhdp
containing $i_0i_1\ldots i_j$; {\bf (3)} once arrived to the first vertex $v$
of an mhdp, or in the case that $u=v$ is such a first vertex, consider its
vertical predecessor, that is the tail $z$ of the vertical arc in
$\Lambda_n$ with head $v$, (which is in the column previous to that
containing $v$); {\bf (4)} set $u=z$ and repeat item (2); {\bf (5)} continue until vertex $0$ is reached.

\example  Let $i_0i_1\ldots i_j=2468aa$ be the vertex of $\Lambda_{11}$
whose weight is 15, (the one left out of the encased table above). This is the first (and only) vertex of its (improper) mhdp.  Its
vertical predecessor, in column 14, is $2468a$.  This is preceded horizontally by
24689 and this by 24688, in respective columns 13 and 12.  The vertical
predecessor of 24688 is 2468, in column 11, preceded horizontally by 2467 and
this by 2466, in respective columns 10 and 9.  The vertical predecessor of 2466
is 246, in column 8, preceded horizontally by 245 and this by 244, in respective
columns 7 and 6.  The vertical predecessor of 244 is 24, in column 5, preceded
horizontally by 23 and this by 22, in respective columns 4 and 3.  The vertical
predecessor of 22 is 2, in column 2, preceded horizontally by 1 and this by
$0=u_0$, in respective columns 1 and 0.  Thus we get the following path, with commas
replaced by superindices $m_u$, for horizontal-arc tails $u$, and subindices
$d_u$, for vertical-arc tails $u$, respectively:

$$_{0^{10}1^92_222^823^724_4244^6245^5246_62466^42467^32468_824688^224689^12468a_{
10}2468aa.}$$ We arrive at $c(2468aa)=9\times 7\times 5\times 3$.

This generalizes to the following statement.

\begin{thm} If $n=2k+1$, then the paths realizing the diameter $D(ST_n)$ of $ST_n$ and
starting at $12\ldots n$ end up at exactly $(n-2)(n-4)\ldots 3$
vertices $u$ of the form $\Sigma(u)=\sigma_1\sigma_2\ldots\sigma_n$, with
$\sigma_1=1$ and $\Pi(u)$ expressible as a product of $k+1$ independent
transpositions.  \end{thm}

A string $i_0i_1\ldots i_j$ is said to be admissible if $u_{i_0i_1\ldots i_j}$
is a vertex of $\Lambda_\infty$. Given a positive integer $\omega\leq D(ST_n)$, we
want first to find an expression for the cardinality of the set $V_\omega$ of
vertices of $\Lambda_\infty$ having $\omega$ as their weight in $ST_{\omega+1}$.
Toward this end, we start exemplifying some sequences of admissible
strings for lower values of $\omega$, where subindex strings $i_0i_1\ldots u_j$ of
vertices $u_{i_0i_1\ldots u_j}$ are expressed in a suitable order without commas and  employing the following shorthand dot-notation rule
for certain subsequences:
let $i_0i_1\ldots i_{k-1}i_k.i_{k+1}\ldots i_{j-1}i_j$ stand for the subsequence
composed by all the admissible strings
$i_0i_1\ldots i_{k-1}\iota_k.\iota_{k+1}\ldots \iota_{j-1}i_j$ in
$\Lambda_n$ with $\iota_\ell\geq i_\ell$, for $k\leq\ell<j$.

\examples Some subsequences of admissible strings in $\Lambda_\infty$ are:
$$\begin{array}{||lll||lll|l||}\hline\hline
^{2.2}_{24.4}  &^{=}_{=}&^{\{22\}}_{\{244\}}  &^{2.i_1}_{24.i_2}&^{=}_{=}&^{\{2i_1,3i_1,\ldots,i_1i_1\}}_{\{24i_2,25i_2,\ldots,2i_2i_2\}}&^{i_1>2}_{i_2>4}\\
^{36.6}_{246.6}&^{=}_{=}&^{\{366\}}_{\{2466\}}&^{36.i_2}_{246.i_3}&^{=}_{=}&^{\{36i_2,37i_2,\ldots,3i_2i_2\}}_{\{246i_3,247i_3,\ldots,24i_3i_3\}}&^{i_2>6}_{i_3>6}\\
^{369.9}       &^{=}    &^{\{3699\}}          &^{369.i_3}&^{=}&^{\{369i_3,36ai_3,\ldots,36i_3i_3\},}&^{i_3>9}\\\hline\hline
\end{array}$$
 For $\omega=0,1,\ldots,15=f$ we can express $V_\omega$ as
follows, where hexadecimal notation is used:
$$\begin{array}{|llllllll|}\hline
^{V_0=\{0\}}_{V_1=\{1\}}&&&&&&^{V_3=\{3,}_{V_4=\{4,}&^{2.2\}}_{2.3\}}\\
^{V_2=\{2\}}&&&&&&^{V_5=\{5,}&^{2.4\}}\\
^{V_6=\{6,}_{V_7=\{7,}&^{2.5,}_{2.6,}&^{24.4\}}_{24.5\}}&&&&&\\
^{V_8=\{8,}_{V_9=\{9,}&^{2.7,}_{2.8,}&^{24.6,}_{24.7,}&^{36.6\}}_{36.7,}&&&&\\
_{V_a=\{a,}&_{2.9,}&^{246.6\}}_{24.8,}&_{36.8,}&_{48.8,}&&&\\
_{V_b=\{b,}&_{2.a,}&^{246.7,}_{24.9,}&^{257.7\}}_{36.9,}&_{48.9,}&&&\\
_{V_c=\{c,}&_{2.b,}&^{246.8,}_{24.a,}&^{257.8,}_{36.a,}&^{268.8\}}_{48.a,}&_{5a.a,}&&\\
&&^{246.9,}_{369.9}&^{257.9,}&^{268.9,}&^{279.9,}&&\\
_{V_d=\{d,}&_{2.c,}&^{2468.8\}}_{24.b,}&_{36.b,}&_{48.b,}&_{5a.b,}&&\\
&&^{246.a,}_{369.a,}&^{257.a,}_{37a.a,}&^{268.a,}&^{279.a,}&^{28a.a,}&\\
_{V_e=\{e,}&_{2.d,}&^{2468.9,}_{24.c,}&^{2579.9\}}_{36.c,}&_{48.c,}&_{5a.c,}&_{6c.c,}&\\
&&^{246.b,}_{369.b,}&^{257.b,}_{37a.b,}&^{268.b,}_{38b.b,}&^{279.b,}&^{28a.b,\,29b.b,}&\\
_{V_f=\{f,}&_{2.e,}&^{2468.a,}_{24.d,}&^{2579.a,}_{36.d,}&^{268a.a\}}_{48.d,}&_{5a.d,}&_{6c.d,}&\\
&&^{246.c,}_{369.c,}&^{257.c,}_{37a.c,}&^{268.c,}_{38b.c,}&^{279.c,}_{39c.c}&^{28a.c,29b.c,\,2ac.c,}&\\
&&^{2468.b,}_{2468a.a\}}&^{2579.b,}&^{268a.b,}&^{279b.b,}&&\\\hline
\end{array}$$
The ten last $V_\omega$ here are expressible as:
$$\begin{array}{|lll|llll|}\hline
^{V_6=\{6,}_{V_7=\{7,}&^{2.5,}_{2.6,}&^{2.44\}}_{2.45\}}&^{V_9=\{9,}_{V_a=\{a,}&^{2.8,}_{2.9,}     &^{2.47,}_{2.48,}&^{2.466\}}_{2.467\}}\\
^{V_8=\{8,}&^{2.7,}&^{2.46\}}  &^{V_b=\{b,}  &^{2.a,}&^{2.49,}        &^{2.468\}}           \\\hline
\end{array}$$\vspace*{-3mm}
$$\begin{array}{|lllllll|}
^{V_c=\{c,}_{V_d=\{d,}           &^{2.b,}_{2.c,}       &^{2.4a,}_{2.4b,}            &^{2.469,}_{2.46a,}             &^{2.4688\}}_{2.4689\}}       &                &                      \\
^{V_e=\{e,}_{V_f=\{f,}           &^{2.d,}_{2.e,}       &^{2.4c,}_{2.4d,}            &^{2.46b,}_{2.46c,}             &^{2.468a\}}_{2.468b,}       &_{2.468aa\}}&                      \\\hline
\end{array}$$
Let $V_\omega^{i_0}$ be the subset of strings of $V_\omega$ starting at $i_0$. We draw the following conclusions, where the dot-notation rule is used:
{\bf(1)} $\lambda=1$ happens in $V_\omega$ just for each subsequence $\omega\geq 0$, and only in $V_\omega^1$; {\bf(2)}
$\lambda=2$ happens in $V_\omega$ for the members of $2.(\omega-1)$, where $\omega\geq 3$, and only in $V_\omega^2$;
{\bf(3)} $\lambda=3$ happens in $V_\omega$:
{\bf(a)} for the members of $24.(\omega-2)$, where $\omega\geq 6$, and only in $V_\omega^2$;
{\bf(b)} for the members of $36.(\omega-2)$, where $\omega\geq 8$, and only in $V_\omega^3$;
$\ldots$
{\bf(z)} for the members of $k(2k).(\omega-2)$, where $\omega\geq 2(k+1)$, and only in $V_\omega^k,$ $(k\geq 2)$;
{\bf(4)} $\lambda=4$ happens in $V_\omega$:
{\bf(a)} for the members of $246.(\omega-3)$ where $\omega\geq 9$, and only in $V_\omega^2$;
{\bf(b)} for the members of $369.(\omega-3)$ where $\omega\geq 12$, and only in $V_\omega^3$;
$\ldots$
{\bf(z)} for the members of $k(2k)(3k).(\omega-3)$, where $\omega\geq 3(k+1)$, and only in $V_\omega^k$, $(k\geq 2)$.
The following result is obtained.

\begin{thm} {\bf(a)} $\lambda=1$ happens in $V_\omega$, and only for the
strings of $V_\omega$ starting at $i_0$; {\bf(b)} for each $k\geq 2$, any fixed
$\lambda >1$ happens in $V_\omega$ for the members of
$k(2k)(3k)\ldots((\lambda-1)k).(\omega-\lambda+1)$,
where $\omega\geq (\lambda-1)(k+1)$, and only for the subsets $V_\omega^k$.  \end{thm}

Let $W_\omega^k\subseteq V_\omega$ consist of the
strings of length $\lambda=k$ in the statement of Theorem 3. Then
$|W_\omega^1|=1 $ and $|W_\omega^k|= 0$ whenever $\omega<3k$, for $k\geq 2$.
Moreover, if
$S_j^0=1$ and $S_j^h=\sum_{k=1}^jS_k^{h-1}$, ($h>0$), for every $j\geq 1$, so $S_j^h-S_{j-1}^h+S_j^{h-1}$ for $h>0$ and $j>1$, then
\begin{equation}S_j^h={j+h-1\choose h},\;\;\;\;\;|W_1^k|=S_1^k=k\;\;\;\mbox{ and in general}\end{equation}
\begin{equation}|W_\omega^k|=\sum_{i=0}^{\lfloor\frac{\omega}{k+1}\rfloor}S_{\omega-i(k+1)}^k =\sum_{i=0}^{\lfloor\frac{\omega}{k+1}\rfloor}{\omega-ik-i+k-1\choose k},\end{equation}
for every weight $\omega$ valid in $ST_{\omega+1}$
and every string length $k$.

\begin{thm} For $0<\omega\in\Z$, the number of vertices of
$ST_{\omega+1}$ having weight $\omega$
is given by the finite sum
$|V_\omega|=|W_\omega^1|+|W_\omega^2|+\ldots +|W_\omega^k|+\ldots.$ \end{thm}

It is easy to establish the following expression for the diameter
$D(n)=D(ST_n)$ of $ST_n$.

\begin{propo}
The diameter of $ST_n$ is $D(n)=\lfloor\frac{n-1}{2}\rfloor+n-1$.
\end{propo}

Let $V_\omega(n)$ be the set of vertices of $\Lambda_n$ having weight
$\omega$. Let $W_\omega^k(n)$ be the subset of admissible strings
corresponding to vertices of $V_\omega(n)$ whose length $\lambda$ is equal to
$k$. Then, from the tables ${\mathcal T}_n$ we get:
\begin{equation}|W_\omega^k(n)|=|W_\omega^k|,\;(0\leq k<n);\end{equation}
\begin{equation}|W_\omega^k(n)|=|W_\omega^k|-\sum_{j=0}^{k-n}|W_\omega^j|,\;(n\leq k\leq D(n)).\end{equation}
The main result of the section follows.

\begin{thm} The cardinality of the set of vertices of $ST_n$ having weight
$\omega$ is
$$|V_\omega(n)|=|W_\omega^0(n)|+|W_\omega^1(n)|+\ldots+|W_\omega^{D(n)}(n)|=\sum_{i=0}^{D(n)}W_\omega^i,$$
where the terms of the displayed sum are obtained by means of equations $(1)$, $(2)$, $(3)$ and $(4)$ presented above.
\end{thm}

\proof The equations and the statement of the theorem arise naturally from the patterns in the tables ${\mathcal T}_n$ and the previous results.  \qfd

\section{Weight distributions of E-sets in $ST_n$}

It was proved in \cite{io} that if $1\leq i\leq n$, then, the vertex subset
$C_i$ of $ST_n$ corresponding to the permutations
$\sigma_1\sigma_2\ldots\sigma_n$ with a fixed $\sigma_1=i$ forms an E-set.
This is the only way of getting an E-set in $ST_n$.
Furthermore, it can be seen that the E-sets of $ST_n$ form a partition of the vertex set of $ST_n$.

Having established in Section 4 the distribution of weights of vertices of $ST_n$, we ask, How does such a distribution restricts to each $C_i$?

\begin{propo} The vertices $u$ of $\Lambda_n$ with
$\Sigma(u)=\sigma_1\sigma_2\ldots\sigma_n$ and
$\sigma_1=1$ represent all the vertices of $ST_n$ with $\sigma_1=1$.
They have associated admissible strings $i_0i_1\ldots
i_{j-1}i_j$ with $i_{j-1}=i_j$.
\end{propo}

\proof This is clear from the developments above. \qfd

Let $V_\omega^i(n)$ be the set of vertices of $C_i$ having weight $\omega$ in
$ST_n$, for $1\leq i\leq n$.

\begin{thm}
The weight distribution of the subsets $C_i$ of $ST_{n+1}$, for
$2\leq i\leq n+1$, is given by:
$$\begin{array}{ll}
|V_0^i(n+1)|=0; \\
|V_\omega^i(n+1)|=|V_{\omega-1}(n)|, &
\mbox{for }\omega=1,2,\ldots,2\lfloor\frac{D(n+1)}{2}\rfloor; \\
|V_{D(n+1)}(n+1)|=0, &
\mbox{for }n\mbox{ even, {\rm(}only case not covered above{\rm)}} \\
\end{array}$$
\end{thm}

\proof For each $i\in\{2,\ldots,n+1\}$, the permutations
$\sigma_1\sigma_2\ldots\sigma_{n+1}$ with $\sigma_i=i$ induce a copy $H_i$ of
$ST_n$ in $ST_{n+1}$ containing the identity permutation $12\ldots(n+1)$.
Each vertex $h$ of $H_i$ has a unique neighbor $h^i$ in
$ST_{n+1}\backslash H_i$. Then the collection of all $h^i$ is $C_i$, for each
$i\in\{2,\ldots n+1\}$ fixed.  \qfd

\rema  According to Theorem 8, the $n$ vertex subsets $C_i$ in
$ST_{n+1}$ with $1<i\leq n+1$ have equivalent weight distributions.  Thus, by
multiplying the quantities obtained in the theorems by $n$ and substracting the
results correspondingly from those obtained for $ST_{n+1}$, the
case for $C_1$ can be obtained, which uses that if $n$ is odd then $|V_{D(n)}|=(n-2)(n-4)\ldots\times 5\times 3$, by Theorem 2.

\section{Threading $\Lambda_n$ into an orientation of $ST_n$}

We now modify the Pruning Algorithm into a threading algorithm in order to produce
an orientation $\Gamma_n$ of $ST_n$ whose vertices are those of $\Lambda_n$
(remaining after applying the algorithm) and whose arc set contains the arc
set of $\Lambda_n$.

The Threading Algorithm consists in running the Pruning
Algorithm (on the previously defined $\Lambda_n$), checking whether the last
field $b_ua_u$ of each line in the table ${\mathcal P}_n$ that is being
generated has $b_u=0$ and $a_u\geq 2$. If this is the case, then a
{\it thread}, meaning a new arc, is added to $\Lambda_n$ from $u$ to a vertex
$\psi(u)$ determined as follows. It happens that the penultimate field $C(u)$
was present in a previous line of ${\mathcal P}_n$ corresponding to the tail $\phi(u)$ of a vertical arc $e(u)$ of $\Lambda_n$ having head $\psi(u)$. Then $\psi(u)$ is the head of $e(u)$.

\example Working with ${\mathcal P}_9$, the threads appearing by means of the
Thre\-ading Algorithm are departing from the vertices $u$ with subindex strings
35, 46, 47, 57, 58, 68, 257, 268, 368, whose values $C(u)$ are respectively
32, 42, 43, 52, 53, 62, 232, 242, 332 and whose fields
$b_ua_u=0a_u$ have $a_u= 2,2,3,2,3,2,2,2,2$, respectively. But
the vertices $\phi(u)$ with respective subindex strings 25, 26, 37, 27, 38,
28, 257, 268, 368, have the same corresponding values $C(u)$, presented in
${\mathcal P}_9$ in nondecreasing order: 23, 24, 34, 25, 35, 26, 223, 224, 233,
so the corresponding 1-ics's are the same in both cases. We obtain
the desired orientation of $ST_9$ by adding a thread
from each one of the eight mentioned vertices respectively into the vertices
$\psi(u)$ whose subindex strings are 255, 266, 377, 277, 388, 288, 2477, 2488,
2588, which are the heads of the respective arcs $e(u)$ (that departed
from the vertices $\phi(u)$ mentioned above).

\begin{thm} Any pair $(u,\phi(u))$ appearing during the running of the
Thre\-ad\-ing
Algorithm has the vertices $u$ and $\phi(u)$ with $C(u)=C(\phi(u))$, where the order
of the elements on each side of the equality is irrelevant.  Thus, in
the running of the Threading Algorithm, each consideration of a vertex $u$ of
$\Lambda_n$ with $C(u)$ equal to the $C(v)$ of a previously considered vertex
$v=\phi(u)$ determines a thread from $u$ onto the corresponding $\psi(u)$.
\end{thm}

\proof The statement follows from the previous discussion and Theorem 1, item
1.  \qfd

\rema The Threading Algorithm insured by Theorem 9 produces an orientation
$\Gamma_n$ of $ST_n$ whose vertices represent the 1-ics's of the
permutations on $n$ elements, that is each vertex of $\Gamma_n$ represents all
the permutations on $n$ elements having a specific 1-ics, and there
is a bijective correspondence between the vertices of $\Lambda_n$ and the
1-ics's of permutations on $n$ elements. Thus $\Gamma_n$ may be
referred to as the 1-ics orientation of $ST_n$. Each arc of
$ST_n$ projects into a specific arc of $\Gamma_n$. We still consider that the
arcs of $\Gamma_n$ are `horizontal' and `vertical', as in the case of
$\Lambda_n$, where threads of $\Gamma_n$ are `vertical'. Moreover, the
vertices and arcs of $\Gamma_n$ may be considered as preserving the indications
they inherit from $\Lambda_n$, including the threads, which preserve the
indications of the arcs removed by the Pruning Algorithm. As said above,
the indications of horizontal arcs are of the form $\bullet\ell_u$, so we still have
that the orientations $\Gamma_n$ form a nested sequence of indicated digraphs
and that their limit indicated digraph $\Gamma_\infty$ is well defined and
constitutes a universal graph for this situation.  This corresponds to the
infinite star graph $ST_\infty$ that can be defined as the Cayley graph of the
symmetric group $S_\infty$ with respect to the set of transpositions
$\Theta_\infty=\{(1\;i),\; i=2,\ldots n,\ldots\}$.

\begin{thm}
$\Gamma_n$ can be interpreted as an orientation of $ST_n$ via the
map $\Phi_n:ST_n\rightarrow\Lambda_n$ given by $\Phi_n^{-1}(u)=$
$\rho$-equivalence class of $\Sigma(u)$, for each vertex $u$ of $\Lambda_n$.
Then: {\bf(1)} the value $c(u)$ of each vertex
$u$ of $\Gamma_n$ is the cardinality of $\Phi_n^{-1}(u)$ and
{\bf(2)} the inverse image $\Phi_n^{-1}$ of an
horizontal, {\rm(}vertical{\rm)}, arc $e$ of $\Lambda_n$ is formed by
$c(u^e)$, $(c(u_e))$, arcs subdivided into $c(u^e)/m_e$, $(c(u_e)/d_e)$,
subsets of $m_e$, $(d_e)$, arcs incident each to a common corresponding
vertex in $\Phi_n^{-1}(u_e)$, $(\Phi_n^{-1}(u^e))$.
\end{thm}

\end{document}